\documentclass[12pt]{article}
\usepackage{amssymb,latexsym}
\usepackage{amsfonts,array}

\input epsf

\def\polyY{{\Upsilon}}
\def\Wgen{{\Theta}}

\def\W2xi{\Wgen_{2,\underline{x}}}

\def\deg{\mathop{\rm deg}\nolimits}

\def\End{\mathop{\rm End}\nolimits}

\def\Hom{\mathop{\rm Hom}\nolimits}
\def\gr{\mathop{\rm gr}\nolimits}

\def\id{\mathop{\rm id}\nolimits}
\def\Gal{\mathop{\rm Gal}\nolimits}

\def\A{{\cal A}}

\def\K{F}

\def\V1{\polyY}

\def\realize{\psi}

\def\edge[#1,#2]{]#1;#2[}

\def\V{{\cal V}}
\def\A{{\cal A}}

\def\deg{{\rm deg}}

\def\spl2{{\rm sl_2}}

\def\id{\rm id}

\def\End{\rm End}

\def\btimes{\begin{picture}(0.55,0.5)(-0.4,0)
\put(-0.4,0){\makebox(0,0)[lb]{\smash{\SetFigFont{12}{14.4}{\familydefault}{\mddefault}{\updefault}$\Box$}}}
\put(-0.4,0.06){\makebox(0,0)[lb]{\smash{\SetFigFont{12}{14.4}{\familydefault}{\mddefault}{\updefault}$\times$}}}
\end{picture}}

\def\rank{{\rm rank}}
\def\Hom{{\rm Hom}}

\def\Br3{{\bf Br}_3}
\def\bino(#1,#2){{\left({#1\atop #2}\right)}}

\def\id{{\rm id}}

\def\mod{{\rm mod}}

\def\V{{\cal V}}

\def\Ab{{\bar{\cal A}}}

\def\Abn1{{\Ab_n}}

\def\Hvn1{{{\cal H}_{n,1}}}
\def\Fvn1{{{\cal F}_{n,1}}}

\def\Hn1{{{\cal H}_n'}}
\def\Kn1{{{\cal F}_n'}}
\def\Fn1{{{\cal F}_n'}}
\def\H61{{{\cal H}_6'}}
\def\K61{{{\cal F}_6'}}

\def\bline[#1,#2]{
\thicklines
\ifnum#1>#2\bline[#2,#1]
\else\ifnum#1<#2
\ifcase#1
\ifcase#2
\or\qbezier[60](-1,-1)(-1,-0.5)(-0.5,-0.5)
\qbezier[60](-0.5,-0.5)(0,-0.5)(0,-1)
\or\qbezier[80](-1,-1),(-1,-0.5)(0,-0.5)
\qbezier[80](0,-0.5)(1,-0.5)(1,-1)
\or\put(-1,-1){\line(1,1){2}}
\or\put(-1,-1){\line(1,2){1}}
\or\put(-1,-1){\line(0,1){2}}
\fi
\or\ifcase#2\or
\or\qbezier[80](0,-1),(0,-0.5)(0.5,-0.5)
\qbezier[80](0.5,-0.5)(1,-0.5)(1,-1)
\or\put(0,-1){\line(1,2){1}}
\or\put(0,-1){\line(0,1){2}}
\or\put(0,-1){\line(-1,2){1}}
\fi
\or\ifcase#2\or\or
\or\put(1,-1){\line(0,1){2}}
\or\put(1,-1){\line(-1,2){1}}
\or\put(1,-1){\line(-1,1){2}}
\fi
\or\ifcase#2\or\or\or
\or\qbezier[60](1,1)(1,0.5)(0.5,0.5)
\qbezier[60](0.5,0.5)(0,0.5)(0,1)
\or\qbezier[80](1,1),(1,0.5)(0,0.5)
\qbezier[80](0,0.5)(-1,0.5)(-1,1)
\fi
\or\ifcase#2\or\or\or\or
\or\qbezier[60](-1,1)(-1,0.5)(-0.5,0.5)
\qbezier[60](-0.5,0.5)(0,0.5)(0,1)
\fi\fi\fi\fi
}

\def\epsfigbracket[#1=#2]{\epsffile{ps/#2}}
\def\epsfig#1{\epsfigbracket[#1]}
\def\includegraphics#1{\epsffile{eps/#1}}

\newcommand{\vcentre}[1]{\setbox1=\hbox{\input{fig/#1}}
$\vcenter{\box1}$}

\def\labelled{{labeled}}  \def\st{{\rm st}}
\def\rank{{\rm rank}} \def\Hom{{\rm Hom}}
\def\sigmab{{\bar{\sigma}}}

 \def\V{{\cal V}}

  \def\sour{{\rm source}}
\def\tar{{\rm target}} \def\na{{\rm na}} \def\X{{X}} \def\V{{\cal
V}}  
\def\H{{\widetilde{H}}}  
\def\A{{\cal A}} \def\Ab{{\bar{\A}}} \def\L{{\cal L}}
 \def\K{{\cal K}} 
 \def\gr{{\rm gr}}

\def\id{{\rm id}}  
 \def\mod{{\rm mod}} \def\realize{\psi}
\def\Ker{{\rm Ker}} \def\End{{\rm End}} \def\Grs{{\rm Gr}^\pm}

 \newtheorem{lemma}{Lemma}
\newtheorem{prop}[lemma]{Proposition}

\newtheorem{theorem}[lemma]{Theorem}
\newtheorem{theoremX}{Theorem} \newtheorem{coro}[lemma]{Corollary}
\newtheorem{defi}[lemma]{Definition}
\newtheorem{fact}[lemma]{Fact}

\newcommand{\N}{{\ensuremath{\mathbb N}}}
\newcommand{\Z}{{\ensuremath{\mathbb Z}}}
\newcommand{\Q}{{\ensuremath{\mathbb Q}}}
\newcommand{\R}{{\ensuremath{\mathbb R}}}

\def\DottedCircle{
\bezier{4}(0.966,-0.259)(1.04,0)(0.966,0.259)
\bezier{4}(0.966,0.259)(0.897,0.518)(0.707,0.707)
\bezier{4}(0.707,0.707)(0.518,0.897)(0.259,0.966)
\bezier{4}(0.259,0.966)(0,1.04)(-0.259,0.966)
\bezier{4}(-0.259,0.966)(-0.518,0.897)(-0.707,0.707)
\bezier{4}(-0.707,0.707)(-0.897,0.518)(-0.966,0.259)
\bezier{4}(-0.966,0.259)(-1.04,0)(-0.966,-0.259)
\bezier{4}(-0.966,-0.259)(-0.897,-0.518)(-0.707,-0.707)
\bezier{4}(-0.707,-0.707)(-0.518,-0.897)(-0.259,-0.966)
\bezier{4}(-0.259,-0.966)(0,-1.04)(0.259,-0.966)
\bezier{4}(0.259,-0.966)(0.518,-0.897)(0.707,-0.707)
\bezier{4}(0.707,-0.707)(0.897,-0.518)(0.966,-0.259)
}
\def\Endpoint[#1]{
}
\def\Arc[#1]{
\thicklines                     
\ifcase#1
\bezier{25}(0.966,-0.259)(1.04,0)(0.966,0.259)
\or
\bezier{25}(0.966,0.259)(0.897,0.518)(0.707,0.707)
\or
\bezier{25}(0.707,0.707)(0.518,0.897)(0.259,0.966)
\or
\bezier{25}(0.259,0.966)(0,1.04)(-0.259,0.966)
\or
\bezier{25}(-0.259,0.966)(-0.518,0.897)(-0.707,0.707)
\or
\bezier{25}(-0.707,0.707)(-0.897,0.518)(-0.966,0.259)
\or
\bezier{25}(-0.966,0.259)(-1.04,0)(-0.966,-0.259)
\or
\bezier{25}(-0.966,-0.259)(-0.897,-0.518)(-0.707,-0.707)
\or
\bezier{25}(-0.707,-0.707)(-0.518,-0.897)(-0.259,-0.966)
\or
\bezier{25}(-0.259,-0.966)(0,-1.04)(0.259,-0.966)
\or
\bezier{25}(0.259,-0.966)(0.518,-0.897)(0.707,-0.707)
\or
\bezier{25}(0.707,-0.707)(0.897,-0.518)(0.966,-0.259)
\fi}
\def\DottedArc[#1]{
\ifcase#1
\bezier{4}(0.966,-0.259)(1.04,0)(0.966,0.259)
\or
\bezier{4}(0.966,0.259)(0.897,0.518)(0.707,0.707)
\or
\bezier{4}(0.707,0.707)(0.518,0.897)(0.259,0.966)
\or
\bezier{4}(0.259,0.966)(0,1.04)(-0.259,0.966)
\or
\bezier{4}(-0.259,0.966)(-0.518,0.897)(-0.707,0.707)
\or
\bezier{4}(-0.707,0.707)(-0.897,0.518)(-0.966,0.259)
\or
\bezier{4}(-0.966,0.259)(-1.04,0)(-0.966,-0.259)
\or
\bezier{4}(-0.966,-0.259)(-0.897,-0.518)(-0.707,-0.707)
\or
\bezier{4}(-0.707,-0.707)(-0.518,-0.897)(-0.259,-0.966)
\or
\bezier{4}(-0.259,-0.966)(0,-1.04)(0.259,-0.966)
\or
\bezier{4}(0.259,-0.966)(0.518,-0.897)(0.707,-0.707)
\or
\bezier{4}(0.707,-0.707)(0.897,-0.518)(0.966,-0.259)
\fi}
\def\Chord[#1,#2]{
\thinlines
\ifnum#1>#2\Chord[#2,#1]
\else\ifnum#1<#2
\ifcase#1
\ifcase#2
\or\qbezier(1,0)(0.516,0.138)(0.866,0.5)
\or\qbezier(1,0)(0.45,0.26)(0.5,0.866)
\or\qbezier(1,0)(0.327,0.327)(0,1)
\or\qbezier(1,0)(0.179,0.311)(-0.5,0.866)
\or\qbezier(1,0)(0.0536,0.2)(-0.866,0.5)
\or\put(1, 0){\line(-2, 0){2}}
\or\qbezier(1,0)(0.0536,-0.2)(-0.866,-0.5)
\or\qbezier(1,0)(0.179,-0.311)(-0.5,-0.866)
\or\qbezier(1,0)(0.327,-0.327)(0,-1)
\or\qbezier(1,0)(0.45,-0.26)(0.5,-0.866)
\or\qbezier(1,0)(0.516,-0.138)(0.866,-0.5)
\fi
\or\ifcase#2\or
\or\qbezier(0.866,0.5)(0.378,0.378)(0.5,0.866)
\or\qbezier(0.866,0.5)(0.26,0.45)(0,1)
\or\qbezier(0.866,0.5)(0.12,0.446)(-0.5,0.866)
\or\qbezier(0.866,0.5)(0,0.359)(-0.866,0.5)
\or\qbezier(0.866,0.5)(-0.0536,0.2)(-1,0)
\or\put(0.866, 0.5){\line(-5, -3){1.73}}
\or\qbezier(0.866,0.5)(0.146,-0.146)(-0.5,-0.866)
\or\qbezier(0.866,0.5)(0.311,-0.179)(0,-1)
\or\qbezier(0.866,0.5)(0.446,-0.12)(0.5,-0.866)
\or\qbezier(0.866,0.5)(0.52,0)(0.866,-0.5)
\fi
\or\ifcase#2\or\or
\or\qbezier(0.5,0.866)(0.138,0.516)(0,1)
\or\qbezier(0.5,0.866)(0,0.52)(-0.5,0.866)
\or\qbezier(0.5,0.866)(-0.12,0.446)(-0.866,0.5)
\or\qbezier(0.5,0.866)(-0.179,0.311)(-1,0)
\or\qbezier(0.5,0.866)(-0.146,0.146)(-0.866,-0.5)
\or\put(0.5, 0.866){\line(-3, -5){1}}
\or\qbezier(0.5,0.866)(0.2,-0.0536)(0,-1)
\or\qbezier(0.5,0.866)(0.359,0)(0.5,-0.866)
\or\qbezier(0.5,0.866)(0.446,0.12)(0.866,-0.5)
\fi
\or\ifcase#2\or\or\or
\or\qbezier(0,1.)(-0.138,0.516)(-0.5,0.866)
\or\qbezier(0,1.)(-0.26,0.45)(-0.866,0.5)
\or\qbezier(0,1.)(-0.327,0.327)(-1,0)
\or\qbezier(0,1.)(-0.311,0.179)(-0.866,-0.5)
\or\qbezier(0,1.)(-0.2,0.0536)(-0.5,-0.866)
\or\put(0, 1){\line(0, -2){2}}
\or\qbezier(0,1.)(0.2,0.0536)(0.5,-0.866)
\or\qbezier(0,1.)(0.311,0.179)(0.866,-0.5)
\fi
\or\ifcase#2\or\or\or\or
\or\qbezier(-0.5,0.866)(-0.378,0.378)(-0.866,0.5)
\or\qbezier(-0.5,0.866)(-0.45,0.26)(-1,0)
\or\qbezier(-0.5,0.866)(-0.446,0.12)(-0.866,-0.5)
\or\qbezier(-0.5,0.866)(-0.359,0)(-0.5,-0.866)
\or\qbezier(-0.5,0.866)(-0.2,-0.0536)(0,-1)
\or\put(-0.5, 0.866){\line(3, -5){1}}
\or\qbezier(-0.5,0.866)(0.146,0.146)(0.866,-0.5)
\fi
\or\ifcase#2\or\or\or\or\or
\or\qbezier(-0.866,0.5)(-0.516,0.138)(-1,0)
\or\qbezier(-0.866,0.5)(-0.52,0)(-0.866,-0.5)
\or\qbezier(-0.866,0.5)(-0.446,-0.12)(-0.5,-0.866)
\or\qbezier(-0.866,0.5)(-0.311,-0.179)(0,-1)
\or\qbezier(-0.866,0.5)(-0.146,-0.146)(0.5,-0.866)
\or\put(-0.866, 0.5){\line(5, -3){1.73}}
\fi
\or\ifcase#2\or\or\or\or\or\or
\or\qbezier(-1,0)(-0.516,-0.138)(-0.866,-0.5)
\or\qbezier(-1,0)(-0.45,-0.26)(-0.5,-0.866)
\or\qbezier(-1,0)(-0.327,-0.327)(0,-1)
\or\qbezier(-1,0)(-0.179,-0.311)(0.5,-0.866)
\or\qbezier(-1,0)(-0.0536,-0.2)(0.866,-0.5)
\fi
\or\ifcase#2\or\or\or\or\or\or\or
\or\qbezier(-0.866,-0.5)(-0.378,-0.378)(-0.5,-0.866)
\or\qbezier(-0.866,-0.5)(-0.26,-0.45)(0,-1)
\or\qbezier(-0.866,-0.5)(-0.12,-0.446)(0.5,-0.866)
\or\qbezier(-0.866,-0.5)(0,-0.359)(0.866,-0.5)
\fi
\or\ifcase#2\or\or\or\or\or\or\or\or
\or\qbezier(-0.5,-0.866)(-0.138,-0.516)(0,-1)
\or\qbezier(-0.5,-0.866)(0,-0.52)(0.5,-0.866)
\or\qbezier(-0.5,-0.866)(0.12,-0.446)(0.866,-0.5)
\fi
\or\ifcase#2\or\or\or\or\or\or\or\or\or
\or\qbezier(0,-1.)(0.138,-0.516)(0.5,-0.866)
\or\qbezier(0,-1.)(0.26,-0.45)(0.866,-0.5)
\fi
\or\ifcase#2\or\or\or\or\or\or\or\or\or\or
\or\qbezier(0.5,-0.866)(0.378,-0.378)(0.866,-0.5)
\fi\fi\fi\fi}
\def\FullChord[#1,#2]{
\Endpoint[#1]
\Endpoint[#2]
\Arc[#1]
\Arc[#2]
\Chord[#1,#2]
}
\def\EndChord[#1,#2]{
\Endpoint[#1]
\Endpoint[#2]
\Chord[#1,#2]
}
\def\Picture#1{
\begin{picture}(2,1)(-1,-0.167)
#1
\end{picture}
}
\def\DottedChordDiagram[#1,#2]{
\Picture{\DottedCircle \FullChord[#1,#2]}
}
\def\ExtChord[#1,#2]{
\Endpoint[#1]\Endpoint[#2]
\thinlines
\ifnum#1>#2\ExtChord[#2,#1]
\else\ifnum#1<#2
\ifcase#1
\or\ifcase#2
\or\or\or\or\or\or\or\or\or

\qbezier[80](0,-1)(-0.1,-1.4)(0.25,-1.35)
\qbezier[80](0.25,-1.35)(1.35,-1.2)(1.35,0)
\qbezier[80](1.35,0)(1.35,0.95)(0.866,0.5) \fi \or\or\ifcase#2
\or\or\or\or\or\or\or\or\or\or\or
\qbezier[80](0,1)(-0.1,1.4)(0.25,1.35)
\qbezier[80](0.25,1.35)(1.35,1.2)(1.35,0)
\qbezier[80](1.35,0)(1.35,-0.95)(0.866,-0.5)
\fi
\or\ifcase#2
\or\or\or\or\or\or\or\or
\qbezier[80](-0.5,0.866)(-0.65,1.1)(-0.75,1.1)
\qbezier[60](-0.75,1.1)(-1.35,1.1)(-1.35,0)
\qbezier[60](-0.5,-0.866)(-0.65,-1.1)(-0.75,-1.1)
\qbezier[80](-0.75,-1.1)(-1.35,-1.1)(-1.35,0)
\fi\fi\fi\fi
}

\begin{document}
\pagestyle{plain}
{\parindent0cm
\setcounter{topnumber}{4}
\setcounter{totalnumber}{6}

\title{\bf Universal Vassiliev invariants of links in coverings of $3$-manifolds}
\author{Jens Lieberum}
\date{}
\maketitle

\begin{abstract}
We study Vassiliev invariants of links in a $3$-manifold~$M$ by
using chord diagrams labeled by elements of the fundamental group
of~$M$. We construct universal Vassiliev invariants of links
in~$M$, where $M=P^2\times [0,1]$ is a cylinder over the real
projective plane~$P^2$, $M=\Sigma\times [0,1]$ is a cylinder over
a surface~$\Sigma$ with boundary, and $M=S^1\times S^2$. A finite
covering~$p:N\longrightarrow M$ induces a map~$\pi_1(p)^*$ between
labeled chord diagrams that corresponds to taking the
preimage~$p^{-1}(L)\subset N$ of a link~$L\subset M$. The
maps~$p^{-1}$ and~$\pi_1(p)^*$ intertwine the constructed
universal Vassiliev invariants.
\end{abstract}

\section*{Introduction}\addcontentsline{toc}{section}
{\numberline{}Introduction}

Let $\Sigma$ be an oriented surface with non-empty boundary and let~$I=[0,1]$.
The universal Vassiliev invariant of links in $I^2\times I$
has been generalized to links in~$\Sigma\times I$
by Andersen, Mattes, and Reshetikhin (\cite{AMR}).
Their universal invariant takes values in a space of chord diagrams,
where by a chord diagram they mean a homotopy class
of certain maps from a usual chord diagram into the
surface~$\Sigma$.
By another approach
a universal Vassiliev invariant of braids in~$\Sigma\times I$ is constructed
in~\cite{GMP},
where~$\Sigma$ is a closed surface of genus~$\geq 1$. This invariant is universal for
Vassiliev invariants with values in an abelian group and separates braids.
Vassiliev invariants of links in a closed oriented $3$-manifold~$M$ are obtained
by using Dehn surgery presentations of~$M$ (\cite{LMO}, \cite{BGRT}).
In this general situation, a universal Vassiliev invariant of links
is not known.

In this paper we start by studying Vassiliev invariants of links
in an arbitrary connected $3$-manifold~$M$. For this purpose we
use a vector space~$\Ab(M)$ spanned by chord diagrams that are
\labelled{} by elements of the fundamental group of~$M$. We
construct a universal Vassiliev invariant of links in~$M$,
where~$M$ is one of the following manifolds:

\begin{itemize}
\item $M=\Sigma\times I$, where~$\Sigma$ is a connected compact
surface with non-empty boundary,

\item $M=P^2\times I$,
where $P^2$ is the real projective plane,

\item $M=S^1\times S^2$.
\end{itemize}

In the first two cases the universal Vassiliev invariant takes
values in a completion of~$\Ab(M)$, whereas in the third case it
takes values in a {\em quotient} of the completion. In the
construction of the universal Vassiliev invariant of links
in~$\Sigma\times I$ we use a decomposition of the surface~$\Sigma$
into an oriented disk and bands that are glued to the boundary of
this disk. In difference to~\cite{AMR} we include the case of
non-orientable surfaces~$\Sigma$. The cases~$M=P^2\times I$ and
$M=S^1\times S^2$ are treated by representing links in~$M$ by
links in~$\Sigma\times I$, where~$\Sigma$ is the M\"obius strip or
$\Sigma=S^1\times I$, respectively.

We extend the definition of the universal Vassiliev invariant of
links in~$\Sigma\times I$ to $I$-bundles over ribbon
graphs~$\Sigma$ which allows a more flexible choice of diagrams of
links in~$\Sigma\times I$.
 For a finite connected covering~$p: N\longrightarrow M$ of $3$-manifolds we define a
map~$\pi_1(p)^*:\Ab(M)\longrightarrow\Ab(N)$ that corresponds to
taking the preimage~$p^{-1}(L)\subset N$ of a link~$L\subset M$.
For coverings of~$\Sigma\times I$, $P^2\times I$, and $S^1\times
S^2$ the maps~$p^{-1}$ and~$\pi_1(p)^*$ intertwine the universal
Vassiliev invariants that are constructed in this paper.

The paper is organized as follows. In Section~\ref{singlink} we
define the Vassiliev filtration on the space ${\cal L}(M)$ of
links in a $3$-manifold~$M$ with associated graded
space~$\gr\L(M)$ and the algebra of Vassiliev invariants~${\cal
V}(M)$ of links in $M$. In Section~\ref{s:lcd} we introduce the
graded coalgebra of $G$-\labelled{} chord diagrams, where~$G$ is a
group together with a homomorphism~$\sigma:G\longrightarrow \{\pm
1\}$. In the case where~$G=\pi_1(M,*)$ and $\sigma$ is the
orientation character of~$M$, we denote this coalgebra
by~$\Ab(M)$. In Section~\ref{s:psi}, we show that for every
connected $3$-manifold~$M$ the coalgebra~$\Ab(M)$ maps
surjectively onto $\gr\L(M)$. In Section~\ref{s:Z} we state our
main results about the relation between the graded
algebra~$\Ab(M)^*$ dual to $\Ab(M)$ and ${\cal V}(M)$:
for~$M=\Sigma\times I$ with $\partial\Sigma\not=\emptyset$ or
$\Sigma=P^2$ we have $\Ab(M)^*\cong {\cal V}(M)$. For $M=S^1\times
S^2$ the space~${\cal V}(M)$ is isomorphic to a subalgebra
of~$\Ab(M)^*$, which we determine. Sections~\ref{s:compcd}
to~\ref{s:s1s2} are devoted to the proofs of these results. In
Section~\ref{s:vari} we discuss variations on the definition of
universal Vassiliev invariants including the case of links in
$I$-bundles over ribbon graphs~$\Sigma$. In Section~\ref{s:cover}
we prove that the universal Vassiliev invariants constructed in
this paper are compatible with finite coverings.
 \section{Links in a $3$-manifold}\label{singlink}

Throughout this paper we work in the piecewise linear category. A
link is an oriented closed $1$-dimensional submanifold of a
$3$-manifold. A knot is a link that has only one connected
component. A singular link may have transversal double points but
no other singularities. Let $\L(M)$ be the vector space freely
generated by isotopy classes of links in~$M$. We associate an
element of~$\L(M)$ to a singular link~$L$ in the following way: we
choose an arbitrary local orientation of~$M$ at each double point
of~$L$ and then desingularize~$L$ to a linear combination of links
without double points by applying the local replacement rule shown
in Figure~\ref{f:desing} to each double point of~$L$.

\begin{figure}[!h]
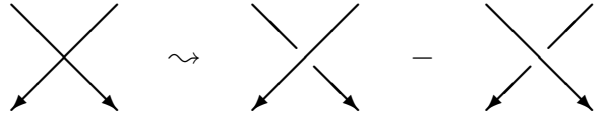

\centering
\Picture{
\thicklines
\put(1,1){\vector(-1,-1){2}}
\put(-1,1){\vector(1,-1){2}}
}
$\quad\leadsto\quad$
\Picture{
\thicklines
\put(1,1){\vector(-1,-1){2}}
\put(-1,1){\line(1,-1){0.83}}
\put(0.17,-0.17){\vector(1,-1){0.83}}
}
$\quad - \quad$
\Picture{
\thicklines
\put(1,1){\line(-1,-1){0.83}}
\put(-0.17,-0.17){\vector(-1,-1){0.83}}
\put(-1,1){\vector(1,-1){2}}
}
\nopagebreak\vspace*{11pt}
\caption{Desingularization of a singular link}\label{f:desing}
\end{figure}

This linear combination is determined by~$L$ and the chosen local
orientations, and it is determined by~$L$ up to a sign. It follows
that we have a well-defined filtration

\begin{equation}
\L(M)=\L_0(M)\supseteq\L_1(M)\supseteq\L_2(M)\ldots
\end{equation}

where~$\L_i(M)$ is spanned by the desingularizations of links in~$M$ with~$i$
double points.
Let $\gr\L(M)=\bigoplus_{n=0}^\infty \L_n(M)/\L_{n+1}(M)$
be the associated graded vector space.
The subspace of~$\gr\L(M)$ generated by
the desingularizations of singular
links with~$\ell$ components and~$n$ double points is denoted
by~$\gr\L(M)_n^\ell$.

The following definitions and results are adapted from~$M=\R^3$
(see \cite{Vo1},\cite{BN1} for details) to the case of an
arbitrary connected~$3$-manifold~$M$. The space~$\L(M)$ has a
coalgebra structure. The comultiplication is given by mapping a
link~$L$ (without double points) to~$L\otimes L$. This coalgebra
structure induces a coalgebra structure on~$\gr\L(M)$.

\begin{defi}
A linear map $v:\L(M)\longrightarrow\Q$ is called a Vassiliev invariant
of degree~$n$ if~$v(\L_{n+1}(M))=0$. Let~$\V_n(M)$ be the vector
space of all Vassiliev invariants of degree~$n$. Define
$\V(M)=\bigcup_{i=0}^\infty \V_i(M)$.
\end{defi}

The spaces $\V_i(M)$ form an increasing sequence

\begin{equation}
\V_0(M)\subseteq\V_1(M)\subseteq\V_2(M)\subseteq\ldots\subseteq\V(M).
\end{equation}

The product $v_1v_2$ of Vassiliev invariants~$v_1$ and~$v_2$ is given by
$(v_1v_2)(L)=v_1(L)v_2(L)$, where $L$ is a link in $M$.
With this product~$\V(M)$ is a subalgebra of~${\cal L}(M)^*$.
If~$v_i$ has degree~$n_i$, then~$v_1v_2$ is a Vassiliev invariant of
degree~$n_1+n_2$.
 \section{Labeled chord diagrams}\label{s:lcd}

Let $G$ be a group. Let $\Gamma$ be a compact one-dimensional
oriented manifold together with a partition of $\partial\Gamma$
into two ordered sets~$s(\Gamma)$ and~$t(\Gamma)$ called {\em
lower} and {\em upper boundary} of $\Gamma$ respectively. For a
set $X$ we denote the set of finite words in the letters~$X$
by~$X^*$. One assigns a symbol~$+$ or~$-$ (resp.~$-$ or~$+$) to a
lower (resp. upper) boundary point of~$\Gamma$ according to
whether~$\Gamma$ directs towards this point or not. Doing this for
all elements of~$s(\Gamma)$ and~$t(\Gamma)$ we obtain two elements
of~$\{+,-\}^*$ called $\sour(\Gamma)$ and $\tar(\Gamma)$
respectively. A {\em $G$-\labelled{} chord diagram}~$D$ with
skeleton~$\Gamma$ consists of the following data:

\begin{itemize}
\item a finite set $S$ of mutually distinct points
on $\Gamma\setminus\partial\Gamma$,

\item a subset $C\subset S$ called {\em chord endpoints},

\item a partition of $C$ into sets of cardinality two called {\em chords},

\item a map from $S\setminus C$ to $G$
that assigns {\em labels} to elements of
  $S\setminus C$.
\end{itemize}

Sometimes we call~$G$-\labelled{} chord diagrams simply chord
diagrams. We consider two chord diagrams~$D$ and~$D'$ with
skeleton~$\Gamma$ and~$\Gamma'$ respectively as being equal, if
there exists a homeomorphism between~$\Gamma$ and~$\Gamma'$ that
preserves all additional data. For a chord diagram with
skeleton~$\Gamma$ we define $\sour(D)=\sour(\Gamma)$ and
$\tar(D)=\tar(\Gamma)$. Define the {\em degree}~$\deg(D)$ of a
chord diagram~$D$ as the number of its chords.

We represent a chord diagram graphically as follows: we represent
the $1$-manifold~$\Gamma$ by drawing the oriented circles and
intervals of $\Gamma$ inside of the strip $\R\times [0,1]$ such
that the lower (resp. upper) boundary points of $\Gamma$ lie on
the horizontal line $\R\times 0$ (resp. $\R\times 1$) and such
that the elements of $s(\Gamma)$ and $t(\Gamma)$ are arranged in
increasing order from left to right. The orientation of the
intervals of~$\Gamma$ is determined by $\sour(\Gamma)$ and
$\tar(\Gamma)$ or indicated by arrows and the circle components
of~$\Gamma$ are oriented counterclockwise in the pictures. The
chords of a chord diagram are represented graphically by
connecting its two endpoints by a line. The $1$-manifold $\Gamma$
is drawn with a thicker pencil than the chords of~$D$. The labels
of a chord diagram are represented by marking the points of
$S\setminus C$ and by writing the labels close to these marked
points. An example of a picture of a $\Z$-\labelled{} chord
diagram~$D$ of degree~$5$ with $\sour(D)=+$ and $\tar(D)=-++$ is
shown in Figure~\ref{f:Dex}.

\begin{figure}[!h]
\centering \setbox1=\hbox{\input{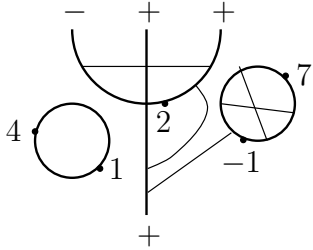}}
$\vcenter{\box1}$
\caption{An example of a picture of a
$\Z$-\labelled{} chord diagram} \label{f:Dex}
\end{figure}

In the following we define vector spaces by generators and
graphical relations. We use the convention that all diagrams in a
graphical relation coincide everywhere except for the parts we
show, and that all configurations of the hidden parts are
possible.

\begin{defi}\label{d:A}
Let $G$ be a group, $\sigma:G\longrightarrow\{\pm 1\}$ a
homomorphism of groups, and let~$s,t\in\{+,-\}^*$. Define
$\A(G,\sigma,s,t)$ as the graded $\Q$-vector space generated by
$G$-\labelled{} chord diagrams~$D$ with~$\sour(D)=s$
and~$\tar(D)=t$ modulo the following relations:

\begin{itemize}
\item The four-term relation (4T):

\nopagebreak

$$
\Picture{
\thicklines
\put(-0.8,1){\vector(0,-1){2}}
\put(0,1){\vector(0,-1){2}}
\put(0.8,1){\vector(0,-1){2}}
\thinlines
\put(-0.8,-0.3){\line(1,0){0.8}}
\put(0,0.4){\line(1,0){0.8}}
}\quad + \quad
\Picture{
\thicklines
\put(-0.8,1){\vector(0,-1){2}}
\put(0,1){\vector(0,-1){2}}
\put(0.8,1){\vector(0,-1){2}}
\thinlines
\put(-0.8,-0.3){\line(1,0){1.6}}
\put(0,0.4){\line(1,0){0.8}}
} \quad = \quad
\Picture{
\thicklines
\put(-0.8,1){\vector(0,-1){2}}
\put(0,1){\vector(0,-1){2}}
\put(0.8,1){\vector(0,-1){2}}
\thinlines
\put(-0.8,0.4){\line(1,0){0.8}}
\put(0,-0.3){\line(1,0){0.8}}
} \quad + \quad
\Picture{
\thicklines
\put(-0.8,1){\vector(0,-1){2}}
\put(0,1){\vector(0,-1){2}}
\put(0.8,1){\vector(0,-1){2}}
\thinlines
\put(-0.8,0.4){\line(1,0){1.6}}
\put(0,-0.3){\line(1,0){0.8}}
}
$$

\bigskip\bigskip

\item The Relation~($\sigma$-Nat):
$\quad\begin{picture}(2,1)(-1,-0)
\thicklines
\put(-1,1){\vector(1,0){2}}
\put(-1,-0.5){\vector(1,0){2}}
\thinlines
\put(-0.4,-0.5){\line(0,1){1.5}}
\put(0.3,-0.5){\circle*{0.15}}
\put(0.3,1){\circle*{0.15}}
\put(0.1,0.5){\makebox(0.4,0.5){$s$}}
\put(0.1,-1){\makebox(0.4,0.5){$s$}}
\end{picture}\ = \ \sigma(s) \
\begin{picture}(2,1.1)(-1,-0.1)
\thicklines
\put(-1,1){\vector(1,0){2}}
\put(-1,-0.5){\vector(1,0){2}}
\thinlines
\put(0.3,-0.5){\line(0,1){1.5}}
\put(-0.4,-0.5){\circle*{0.15}}
\put(-0.4,1){\circle*{0.15}}
\put(-0.6,0.5){\makebox(0.4,0.5){$s$}}
\put(-0.6,-1){\makebox(0.4,0.5){$s$}}
\end{picture}$

\bigskip

\item The Relations~(Rep):
$\quad\begin{picture}(2,0.5)(-1,-0)
\thicklines
\put(-1,0.25){\vector(1,0){2}}
\put(-0.5,0.25){\circle*{0.15}}
\put(0.4,0.25){\circle*{0.15}}
\put(-0.6,-0.3){\hbox{\hss$a$\hss}}
\put(0.3,-0.3){\hbox{\hss$b$\hss}}
\end{picture}\ = \
\begin{picture}(2,0.5)(-1,0)
\thicklines
\put(-1,0.25){\vector(1,0){2}}
\put(0,0.25){\circle*{0.15}}
\put(-0.25,-0.3){\hbox{\hss$ab$\hss}}
\end{picture}\quad , \quad
\begin{picture}(2,0.5)(-1,0)
\thicklines
\put(-1,0.25){\vector(1,0){2}}
\put(0,0.25){\circle*{0.15}}
\put(-0.4,-0.45){\makebox(0.8,0.5){$e$}}
\end{picture}\ = \
\begin{picture}(2,0.5)(-1,0)
\thicklines
\put(-1,0.25){\vector(1,0){2}}
\end{picture}\ ,$

\medskip

where~$e$ is the neutral element of~$G$ and $ab$ denotes the product
of~$a$ and~$b$ in the group~$G$.
\end{itemize}

Let $\Ab(G,\sigma,s,t)$ be the quotient of $\A(G,\sigma,s,t)$ by the framing
independence relation~(FI):
$
\begin{picture}(2,0.5)(-1,0.4)
\thicklines
\put(-1,0.25){\vector(1,0){2}}
\thinlines
\qbezier[70](-0.6,0.25)(-0.6,0.75)(-0.1,0.75)
\qbezier[70](0.4,0.25)(0.4,0.75)(-0.1,0.75)
\end{picture} \quad = \quad 0.
$
\end{defi}

Define the composition $D_1\circ D_2$ of two chord diagrams
with~$\sour(D_1)=\tar(D_2)$ graphically by placing~$D_1$ onto the
top of~$D_2$ and by sticking them together. For a fixed
pair~$(G,\sigma)$ this operation is compatible with all relations
of Definition~\ref{d:A}. For a $G$-\labelled{} chord diagram~$D$
of degree~$n>0$ define

\begin{equation}\label{e:comult}
\Delta(D)=\sum_{D=D'\cup D''} D'\otimes D'',
\end{equation}

where the sum is taken over the~$2^n$ pairs of diagrams~$(D',
D'')$ such that $D'$ and $D''$ have the same \labelled{} skeleton
as~$D$ and the chords of~$D$ are the disjoint union of the chords
of~$D'$ and the chords of~$D''$. For~$D$ with~$\deg(D)=0$
define~$\Delta(D)=D\otimes D$. As in the case $G=\{e\}$ the
spaces~$\A(G,\sigma,s,t)$ and $\Ab(G,\sigma,s,t)$ have a coalgebra
structure with comultiplication~$\Delta$.

Define the category~$\Grs$ as follows.
The objects of~$\Grs$ are
pairs~$(G,\sigma)$ of a group~$G$ and a
homomorphism~$\sigma:G\longrightarrow \{\pm 1\}$.
A morphism~$\varphi:(G,\sigma)
\longrightarrow (H,\sigma')$ in~$\Grs$ is a homomorphism of groups
$\varphi:G\longrightarrow H$ such that $\sigma'=\sigma\circ\varphi$.

Let $s,t\in\{+,-\}^*$ be given. We define a functor
$\bar{\cal F}$ from $\Grs$ to the category of coalgebras
for objects by
$\bar{\cal F}(G,\sigma)=\Ab(G,\sigma,s,t)$. For
morphisms~$\varphi$ we define $\bar{\cal F}(\varphi)=\varphi_*$ by
applying~$\varphi$ to the labels of a chord diagram.
 \section{Mapping chord diagrams to singular links}\label{s:psi}

Denote the empty word by $\emptyset$. We denote
$\Ab(G,\sigma,\emptyset,\emptyset)$ simply by $\Ab(G,\sigma)$. Let
$M$ be a $3$-manifold and $*\in M$. Define a map

\begin{equation}\label{e:sigma}
\sigma:\pi_1(M,*)\longrightarrow\{\pm 1\}
\end{equation}

as follows. Choose a local orientation of~$M$ at the point~$*$. If
the local orientation stays the same when we push it along a
generic representative of an element~$w\in\pi_1(M,*)$, then define
$\sigma(w)=1$, and define $\sigma(w)=-1$ otherwise. The
map~$\sigma$ is called the orientation character of~$M$. We denote
$\Ab(\pi_1(M,*),\sigma)$ by $\Ab(M)$. In this section we
relate~$\Ab(M)$ with~$\gr\L(M)$.

Choose an oriented neighborhood~$U$ of~$*$ with~$U\cong \R^3$. We
identify homotopy classes of paths $w$ from $p_1\in U$ to~$p_2\in
U$ with elements in~$\pi_1(M,*)$ by pulling~$p_1$ and~$p_2$ to~$*$
along paths inside of~$U$. We will use this identification
throughout the rest of this paper. If~$D$ is a
$\pi_1(M,*)$-\labelled{} chord diagram with
$\sour(D)=\tar(D)=\emptyset$, then it is easy to see that one can
construct a singular link~$L_D\subset M$ with the following
properties:

\smallskip

1) All double points of $L_D$ lie inside of $U$.

2) There is an immersion of the skeleton of~$D$ onto~$L_D$ that
fails to be injective only at the preimages of the double points
of~$L_D$. The preimages of double points of $L_D$ correspond to
the pairs of points of~$D$ connected by a chord.

3) Let $w\in \pi_1(M,*)$ be the product of the markings of $D$
along a part of~$D$ between two consecutive (not necessarily
different) chord endpoints. The empty product is considered as
neutral element. By 1) and 2) the segment of~$D$ between these two
chord endpoints is mapped to a part $p$ of $L_D$ connecting two
double points in $U$. We require that the path~$p$ represents
$w\in\pi_1(M,*)$.

4) Let $K$ be a component of~$D$ without chord endpoints. Let
$A\subseteq\pi_1(M,*)$ be the conjugacy class represented by the product of the
labels on~$K$ ($A$ is independent of the choice of the starting point for
this product). Then $K$ is mapped to an unbased loop in~$M$
representing~$A$.

\smallskip

Let us consider an example. Let $M=S^1\times I^2$ and identify
$\pi_1(M,*)$ with~$\Z$. In Figure~\ref{f:bspld} we see a
$\Z$-\labelled{} chord diagram~$D$ and a possible choice of the
singular link~$L_D\subset M$. The neighborhood~$U$ of~$*$ consists
of the right half of~$S^1\times I^2$ in Figure~\ref{f:bspld}.

\begin{figure}[!h]
\centering \setbox1=\hbox{\input{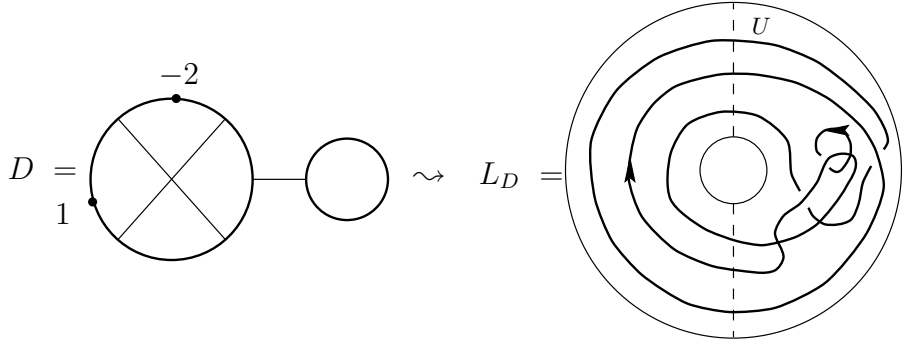}}
$\vcenter{\box1}$ \caption{Mapping a chord diagram $D$ to
a singular link $L_D$}\label{f:bspld}
\end{figure}

Notice that the singular link $L_D$ is not uniquely determined
by~$D$, but we call~$L_D$ every singular link with properties
1)--4). By property 1) we may desingularize all double points of
$L_D$ with the local orientations at double points determined by
the orientation of $U$. We call this desingularization~$[L_D]$.
If~$D$ has~$n$ chords, then we have $[L_D]\in \L_n(M)$.

\begin{prop}\label{diagsing}
Mapping a $\pi_1(M,*)$-\labelled{} chord diagram~$D$ of degree~$n$
to

$$ \psi_M(D)=[L_D]\,\mod\,\L_{n+1}(M) $$

induces a graded surjective morphism of coalgebras
$
{\psi_M}:\Ab(M)\longrightarrow \gr\L(M)$.

\end{prop}
{\bf Proof:} (a) ${\psi_M}$ is well-defined: Consider two singular
links $L$ and $L'$ with the properties 1)--4) of $L_D$. Then one
can pass from $L$ to $L'$ by crossing changes and isotopies. In
other words $[L]-[L']\in\L_{n+1}(M)$. So ${\psi_M}$ is
well-defined for \labelled{} chord diagrams. We have to check that
${\psi_M}$ respects the defining relations of~$\Ab(M)$.

The map ${\psi_M}$ respects the Relations $(4T)$ and $(FI)$
because these relations are a consequence of the Reidemeister
moves in $\R^3\cong U$ and properties 1) and 2) of~$L_D$
(compare~\cite{BN1}). It follows from property 3) of~$L_D$ that
${\psi_M}$ respects the Relations~$(Rep)$.

Call the part of a chord diagram from the left side of the
Relation ($\sigma$--Nat) $P$ and the part from the right side
$P'$. Let~$D$ and~$D'$ be chord diagrams with~$P\subset D$,
$P'\subset D'$ and $D\setminus P=D'\setminus P'$. Let $\hat{s}$ be
an annulus $S^1\times [-1,1]$ in $M$ whose core $S^1\times 0$
represents $s\in\pi_1(M,*)$. We may choose~$L_D$ such that the
oriented arcs of~$P\subset D$ are mapped entirely into~$\hat{s}$
and no other part of~$L_D$ intersects $\hat{s}$. Now we push the
double point of~$L_D$ corresponding to the chord in~$P$
along~$\hat{s}$ until we obtain a singular link~$L_{D'}$. Since
the desingularization of singular links is defined with respect to
the orientation of~$U$, we have
${\psi_M}(D)=\sigma(s){\psi_M}(D')$.

(b) ${\psi_M}$ is surjective: Let $L$ be a singular link with $n$ double
points. Push all double points of $L$ into the neighborhood $U$ of $*$ by an
isotopy and call the resulting singular link $\widetilde{L}$. Let
$D$ be a chord diagram such that~$\widetilde{L}$ has all
properties of $L_D$ (If we replace labels of $D$
between consecutive chord endpoints
by unique labels using the Relation (Rep) and if we fix a choice of
representatives of conjugacy classes as labels of components
of~$D$ without chord endpoints, then the diagram~$D$ is
uniquely determined by $\widetilde{L}$.).
Then~${\psi_M}(D)$ is equal to~$\pm 1$ times the desingularization of~$L$
with freely chosen local orientations at double points modulo
$\L_{n+1}(M)$.

(c) It follows as in the case of $M=\R^3$ that $\psi_M$ is a morphism
of coalgebras (see \cite{Vo1}, Proposition~4).
$\Box$

\medskip

As an immediate application of Proposition~\ref{diagsing} we
obtain the following corollary.

\begin{coro}\label{c:grLM}
(1) If $\pi_1(M,*)$ is finite, then $\gr\L(M)^\ell_n$ is finite
dimensional.

(2) If there exists an element~$s$ with $\sigma(s)=-1$ in the
center of~$\pi_1(M,*)$, then we have~$\L_{2n-1}(M)=\L_{2n}(M)$ for
all~$n>0$.
\end{coro}
{\bf Proof:} (1) If $\pi_1(M,*)$ is finite, then there exists only
a finite number of~$\pi_1(M,*)$-\labelled{} chord diagrams of
degree~$n$ with~$\ell$ circles as skeleton and without neighbored
labels that would allow an application of Relation~$(Rep)$. It
follows from Proposition~\ref{diagsing} that~$\gr\L(M)_n^\ell$ is
finite dimensional in this case.

\smallskip

(2) Let $D$ be a $\pi_1(M,*)$-\labelled{} chord diagram. Let~$s$
be in the center of $\pi_1(M,*)$ with $\sigma(s)=-1$. Using the
Relations~(Rep) of $\Ab(M)$ we insert the labels $s$ and $s^{-1}$
between each pair of consecutive chord endpoints. Then we commute
these labels with the other labels of $D$ such that near the
endpoints of each chord we have labels as shown in
Figure~\ref{f:orirev}.

\begin{figure}[!h]
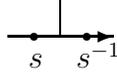

$$ \Picture{ \thicklines \put(-1,0){\vector(1,0){2}} \thinlines
\put(0,0){\line(0,1){0.75}} \put(-0.5,0){\circle*{0.15}}
\put(0.5,0){\circle*{0.15}} \put(-0.6,-0.6){\hbox{\hss$s$\hss}}
\put(0.3,-0.6){\hbox{\hss$s^{-1}$\hss}} } $$ \caption{One can use
the Relation~($\sigma$-Nat) for each chord}\label{f:orirev}
\end{figure}

We apply the Relation~($\sigma$-Nat) for every chord. Using~(Rep)
we forget the pairs $(s^{-1}, s)$ again. This procedure implies
that in $\Ab(M)$ we have $D=(-1)^{\deg D} D$. Hence all elements
of odd degree in~$\Ab(M)$ are~$0$. The surjectivity of ${\psi_M}$
stated in Proposition~\ref{diagsing} implies part~(2) of
Corollary~\ref{c:grLM}. $\Box$

\medskip

The reduction of the coefficient of~$z^{2n+1}$ modulo~$2$ of the
Conway polynomial of links in~$P^2\times I$ (see~\cite{Li2}) is an
example of a Vassiliev invariant of degree~$2n+1$ that is not a
Vassiliev invariant of degree~$2n$. Corollary~\ref{c:grLM} implies
that this invariant cannot be lifted to a~$\Z$-valued Vassiliev
invariant.

Now let us clarify how ${\psi_M}$ depends on the choices we made.
Let $p_1,p_2\in M$. Choose oriented
neighborhoods $U_i\cong \R^3$ of $p_i$. Then
two maps ${\psi_{M,i}}:\Ab(M,p_i)\longrightarrow \gr\L(M)$ are
defined.
If $M$ is orientable and if $U_1$ and $U_2$ induce different orientations
of $M$, then let $\epsilon$ be~$-1$. Let $\epsilon$ be~$1$ otherwise.

\begin{prop}\label{depchoice}
There exists a path $w$ from $p_2$ to $p_1$, such that the following
diagram is commutative
for~$\varphi:\pi_1(M,p_1)\longrightarrow\pi_1(M,p_2)$
defined by~$\varphi(u)=wuw^{-1}$.

$$
\begin{picture}(11,4.5)(-5.5,-2)
\put(-5.5,1){\makebox(2,1){$\Ab(M,p_1)$}}
\put(3.5,1){\makebox(2,1){$\Ab(M,p_2)$}}
\put(-3,1.4){\vector(1,0){6}}
\put(-4,1){\vector(1,-1){2.5}}
\put(4,1){\vector(-1,-1){2.5}}
\put(-1,-2){\makebox(2,1){$\gr\L(M)$}}
\put(-3.9,-0.75){\makebox(1,1){${\psi_{M,1}}$}}
\put(3,-0.75){\makebox(1,1){${\psi_{M,2}}$}}
\put(-3,1.4){\makebox(6,1){$D\mapsto\epsilon^{\deg D}\varphi_*(D)$}}
\end{picture}
$$
\end{prop}

The proof of the proposition is easy and is left to the reader.

 \section{Main results: universal Vassiliev invariants}\label{s:Z}

Consider the completion $\widehat{\A}(M)=\prod_{i=0}^\infty
\Ab(M)_i$ of~$\Ab(M)$. We write elements of~$\widehat{\A}(M)$ as
formal series. The comultiplication~$\Delta$ of~$\Ab(M)$ extends
to a map

$$
\widehat{\Delta}
:\widehat{\A}(M)\longrightarrow
\widehat{\A}(M)\widehat{\otimes}\widehat{\A}(M)=
\prod_{i=0}^\infty \bigoplus_{j=0}^i
\Ab(M)_i\otimes\Ab(M)_{j-i}.
$$

In Section~\ref{s:psi} we saw that there exists a surjective
map~$\psi_M:\Ab(M)\longrightarrow\gr\L(M)$ for every connected manifold~$M$.
In this section we state
stronger results for particular manifolds~$M$.
The following lemma generalizes the invariant of~\cite{AMR} to
non-orientable surfaces~$\Sigma$.

\begin{lemma}\label{l:Kontsevi}
Let $\Sigma$ be a compact connected surface with non-empty
boundary or let $\Sigma$ be the real projective plane. Then there
exists a linear map

$$
Z_{\Sigma\times I}:\L(\Sigma\times I)\longrightarrow \widehat{\A}
(\Sigma\times I)
$$

such that for every $\pi_1(\Sigma\times I,*)$-\labelled{} chord diagram
$D$ of degree~$n$ with $\sour(D)=\tar(D)=\emptyset$ we have

\begin{equation}\label{e:univ}
Z_{\Sigma\times I}\left([L_D]\right) =
D+\mbox{terms of degree $>n$,}
\label{univ}
\end{equation}

and for every link $L$ in $\Sigma\times I$ we have

\begin{equation}\label{e:groupl}
\widehat{\Delta}\left(Z_{\Sigma\times I}(L)\right)
= Z_{\Sigma\times I}(L)
\widehat{\otimes} Z_{\Sigma\times I}(L)\label{groupl}.
\end{equation}
\end{lemma}

The lemma shall be proven in Sections~\ref{s:ZS} and~\ref{s:rp2}.

Let~$\Ab(M)^*$ be the graded dual algebra of the coalgebra~$\Ab(M)$.
The map~$Z_{\Sigma\times I}$ is called {\em universal Vassiliev
invariant} for a reason given by the second
statement of the following theorem.

\begin{theoremX}\label{th:Kontsevi}
Let $\Sigma$ be a compact surface with non-empty boundary
or let $\Sigma$ be the real projective plane. Then
the map
$\psi_{\Sigma\times I}:\Ab(\Sigma\times I)\longrightarrow
\gr\L(\Sigma\times I)$
is an isomorphism of coalgebras.
The formula $(Z_{\Sigma\times I}^*(w))(L)=w(Z_{\Sigma\times
I}(L))$ defines an isomorphism of algebras

$$
Z_{\Sigma\times I}^*:\Ab(\Sigma\times I)^*\longrightarrow \V(\Sigma\times I).
$$
\end{theoremX}
{\bf Proof:}
In view of Proposition~\ref{diagsing}, it suffices to show that
$\psi_{\Sigma\times I}$ is injective for the first statement.
Equation~(\ref{univ}) implies
that~$Z_{\Sigma\times I}$ induces a map from~$\gr\L(\Sigma\times I)$ to
$\Ab(\Sigma\times I)$ which is the left inverse of
$\psi_{\Sigma\times I}$.
This implies that $\psi_{\Sigma\times I}$ is injective.

\smallskip

If~$w\in\Ab(\Sigma\times I)^*$ has degree~$n$, then
equation~(\ref{univ}) implies that $Z_{\Sigma\times I}^*(w)$ is a
Vassiliev invariant of degree~$n$. If~$v$ is a Vassiliev invariant
of degree~$n$, then the map~$W(v)$ defined for chord diagrams~$D$
of degree~$n$ by~$W(v)(D)=v([L_D])$ is in~$\Ab(\Sigma\times
I)^*_n$. Equation~(\ref{univ}) implies that $w=W(Z^*_{\Sigma\times
I}(w))$ and that~$v-Z^*_{\Sigma\times I}(W(v))$ is a Vassiliev
invariant of degree~$n-1$. This implies the injectivity
of~$Z_{\Sigma\times I}^*$ and its surjectivity by induction.
Equation~(\ref{groupl}) implies that $Z_{\Sigma\times I}^*$ is an
isomorphism of algebras. $\Box$

\medskip

Let~$s$ be an element of a group~$G$.
In Figure~\ref{f:sbull} we introduce a notation for a
part of a $G$-\labelled{} chord diagram that
depends on the orientations of strands with an unspecified orientation
in the figure.

\begin{figure}[!h]
\centering
\setbox1=\hbox{\input{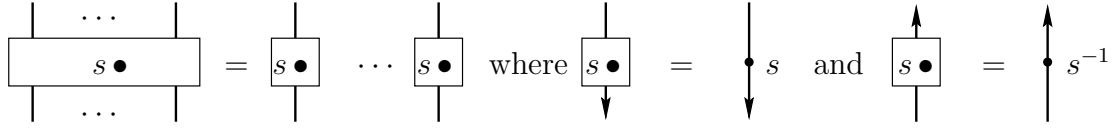}}
$\vcenter{\box1}$
\caption{A notation for oriented \labelled{} strands}\label{f:sbull}
\end{figure}

Similar to Figure~\ref{f:sbull}, we introduce the notation shown in
Figure~\ref{multichord}. The sum in the figure has one term for each
strand of the shown part of a chord diagram.
The free chord end can lead to any point on the chord diagram but this
point must be the same for all terms in the sum.

\begin{figure}[!h]
\centering
\setbox1=\hbox{\input{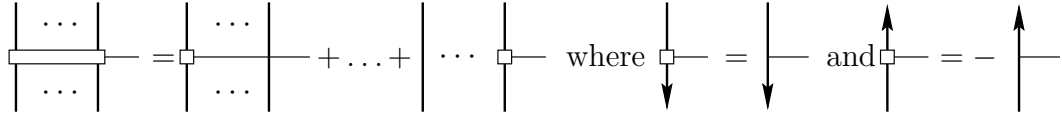}}
$\vcenter{\box1}$
\caption{A notation for certain sums of chord diagrams}\label{multichord}
\end{figure}

We also use the notation obtained by reflection of the diagrams in
Figure~\ref{multichord} in a vertical axis and by combinations of the
notations as shown in the example of Figure~\ref{mchordex}.

\begin{figure}[!h]
\centering
\setbox1=\hbox{\input{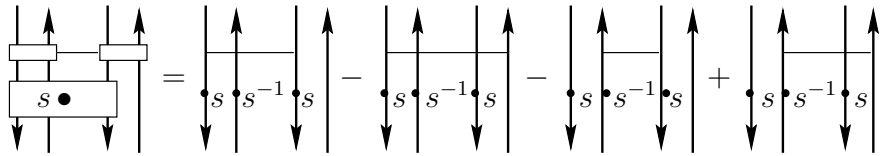}}
$\vcenter{\box1}$
\caption{An example for the notation introduced in
Figures~\ref{f:sbull} and~\ref{multichord}}\label{mchordex}
\end{figure}

Now let~$s$ be a generator of $\pi_1(S^1\times S^2,*)\cong
(\Z,+)$. Let $m$ be the number of strands in
Figure~\ref{mtermrel}. As usual, we assume that the~$m$ terms of
the sum of elements of $\Ab(S^1\times S^2)$ in this figure
coincide except for the shown parts. We further assume that the
hidden part of the diagram in Figure~\ref{mtermrel} has {\em no}
labels\footnote{Lemma~\ref{implrels} implies that
Relation~($S^2$-slide) from Figure~\ref{mtermrel} does not depend
on the choice of~$s$.}.

\begin{figure}[!h]
\centering
\setbox1=\hbox{\input{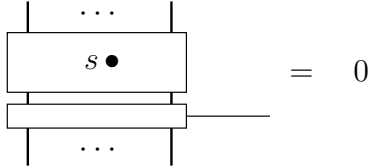}}
$\vcenter{\box1}$
\caption{The Relation~($S^2$-slide)}\label{mtermrel}
\end{figure}

Let ${\cal E}(S^1\times S^2)$ be the quotient of $\Ab(S^1\times
S^2)$ by the Relation~($S^2$-slide) shown in
Figure~\ref{mtermrel}. The canonical projection $p:\Ab(S^1\times
S^2)\longrightarrow {\cal E}(S^1\times S^2)$ has a non-trivial
kernel. For example, the diagram~$D$ in Figure~\ref{f:elex} is
equal to~$0$ in ${\cal E}(S^1\times S^2)$ and it is easy to see
that~$D\not=0$ in $\Ab(S^1\times S^2)$.

\begin{figure}[!h]
\centering
\setbox1=\hbox{\input{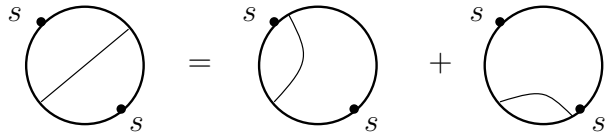}}
$\vcenter{\box1}$
\caption{A non-trivial element in $\Ker(p)$}\label{f:elex}
\end{figure}

In Section~\ref{s:s1s2}
we will construct a map $Z_{S^1\times S^2}$ from
$\L(S^1\times S^2)$ to the completion $\hat{{\cal E}}(S^1\times S^2)$
of~${\cal E}(S^1\times S^2)$
with properties as in
Lemma~\ref{l:Kontsevi}.
This will imply the following theorem.

\begin{theoremX}\label{t:Kon2}
The map $\psi_{S^1\times S^2}$ factors through an isomorphism of
coalgebras

$$
\tilde{\psi}_{S^1\times S^2}:{\cal E}\left(S^1\times S^2\right)\longrightarrow
\gr\L\left(S^1\times S^2\right).
$$

The map $Z_{S^1\times S^2}$ induces an isomorphism of algebras

$$
Z_{S^1\times S^2}^*:{\cal E}\left(S^1\times
  S^2\right)^*\longrightarrow
\V\left(S^1\times S^2\right).
$$
\end{theoremX}

It was shown in~\cite{Eis} that there exist knots in~$S^2\times
S^1$ that can be distinguished by~$\Z/2$-valued Vassiliev
invariants but that cannot be distinguished by $\Q$-valued
Vassiliev invariants.

In Section~\ref{s:cover} we state and prove Theorems~\ref{t:cover}
and~\ref{t:cover2}. By Theorem~\ref{t:cover} a suitably defined
universal Vassiliev invariant of links in $I$-bundles over
surfaces with boundary is compatible with finite coverings. In
Theorem~\ref{t:cover2} the invariant~$Z_{P^2\times I}$ is related
to the usual universal Vassiliev invariant of links in~$\R^3$.
Theorem~\ref{t:cover2} is used in~\cite{Li2} to prove a refinement
of a theorem of Hartley and Kawauchi about the Alexander
polynomial of strongly positive amphicheiral knots in~$\R^3$.

Let us overview the rest of this paper: in Section~\ref{s:compcd}
we derive some general consequences of the defining relations
of~$\A(G,\sigma)$ and some individual properties of~$\Ab(P^2\times
I)$. We introduce a category of tangles in thickened decomposed
surfaces in Section~\ref{s:tangles}. In Section~\ref{s:recallZ} we
recall the definition and some properties of a functor~$Z$ from a
category of non-associative tangles in~$I^2\times I$ to chord
diagrams. In Section~\ref{s:pd} we represent a tangle
in~$\Sigma\times I$ as the composition of a standard tangle
in~$\Sigma\times I$ with a tangle in~$I^2\times I$. By defining
non-associative standard tangles in~$\Sigma\times I$ and by
choosing \labelled{} chord diagrams as their values, we extend the
functor~$Z$ to a functor~$Z_{\cal S}$ from non-associative tangles
in~$\Sigma\times I$ to \labelled{} chord diagrams in
Section~\ref{s:ZS}. This will imply
Lemma~\ref{l:Kontsevi} for surfaces with non-empty boundary. In
Section~\ref{s:rp2} we prove Lemma~\ref{l:Kontsevi}
for~$\Sigma=P^2$. Section~\ref{s:s1s2} is devoted to the proof of
Theorem~\ref{t:Kon2}. The main purpose of Section~\ref{s:vari} is
to generalize the invariant~$Z_{\Sigma\times I}$ to links in
$I$-bundles over ribbon graphs. In Section~\ref{s:cover} we use
this generalization to prove Theorems~\ref{t:cover}
and~\ref{t:cover2} about universal Vassiliev invariants and
coverings.

\section{Computations with chord diagrams}\label{s:compcd}

For later use we provide some combinatorial identities for labeled
chord diagrams in this section. Throughout this section $G$ is a
group and~$\sigma:G\longrightarrow\{\pm 1\}$ is a homomorphism.
The following lemma is a consequence of the four-term relation.

\begin{lemma}\label{implrels}
The relation shown in Figure~\ref{2mtermrel} holds
in~$\A(G,\sigma)$, where $g_1,\ldots,g_m\in G$ are all labels of
the diagrams.

\begin{figure}[!h]
\centering \setbox1=\hbox{\input{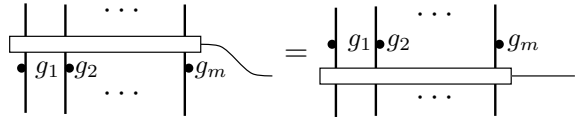}}
$\vcenter{\box1}$ \caption{A relation for chord diagrams
without boundary}\label{2mtermrel}
\end{figure}
\end{lemma}
{\bf Proof:} Let $D$ be a $G$-\labelled{} chord diagram with
$\sour(D)=\tar(D)=\emptyset$ that coincides with the diagrams in
Figure~\ref{2mtermrel} except for the chord shown in this figure.
Keep in mind the common endpoint~$c\in D$ of the chord in
Figure~\ref{2mtermrel}. Consider a picture of~$D$ that is
decomposed by horizontal lines into parts as in Figure~\ref{parts}
(if the label~$g$ is equal to~$e$, we may omit it). Let the
non-trivial labels be the maxima of the picture.

\begin{figure}[!h]
\centering
\setbox1=\hbox{\input{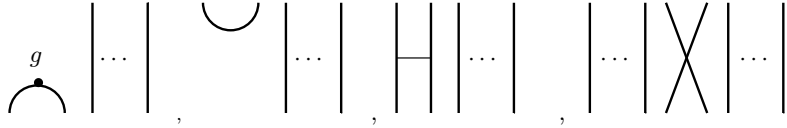}}
$\vcenter{\box1}$
\caption{Parts of a \labelled{} chord diagram}\label{parts}
\end{figure}

We obviously have the relations shown in Figure~\ref{rels}. The
second relation is implied by the $(4T)$-relation in
$\A(G,\sigma)$. \nopagebreak{}
\begin{figure}[!h]
\centering
\setbox1=\hbox{\input{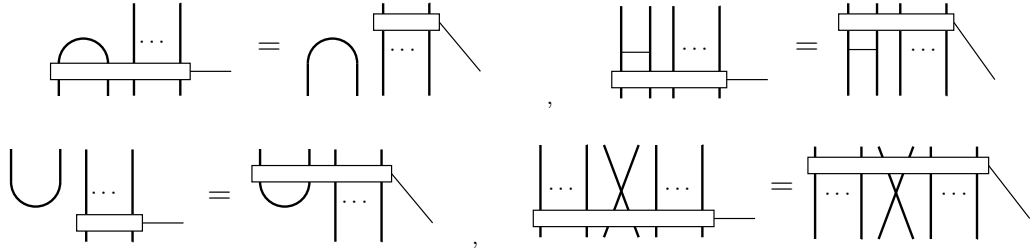}}
$\vcenter{\box1}$
\caption{Relations for parts of a \labelled{} chord diagram}\label{rels}
\end{figure}

Starting with the empty sum at the bottom of the picture of~$D$
and applying the relations from Figure~\ref{rels} for chords
leading to~$c\in D$, we arrive near the top of the diagram as
shown in Figure~\ref{2mtermrelb}. \nopagebreak{}
\begin{figure}[!h]
\centering
\setbox1=\hbox{\input{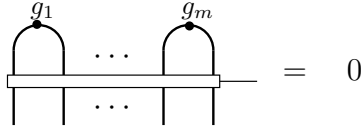}}
$\vcenter{\box1}$
\caption{A consequence of the relations from
Figure~\ref{rels}}\label{2mtermrelb}
\end{figure}

This relation is equivalent to the relation from the lemma.
$\Box$

\medskip

The proof of the next lemma is easy.
It follows from the Relations~(Rep)
and~($\sigma$-Nat) along the same lines as the proof of
Lemma~\ref{implrels}.

\begin{lemma}\label{sDNat}
Let $D$ be a chord diagram without labels. Let $s\in G$.
Define~$\sigmab(s)\in\{0,1\}$ by~$\sigma(s)=(-1)^{\sigmab(s)}$.
Then the relation in Figure~\ref{f:sDNat} holds
in~$\A(G,\sigma,a,b)$, where~$a=\sour(D)$ and $b=\tar(D)$.
\nopagebreak{}
\begin{figure}[!h]
\centering
\setbox1=\hbox{\input{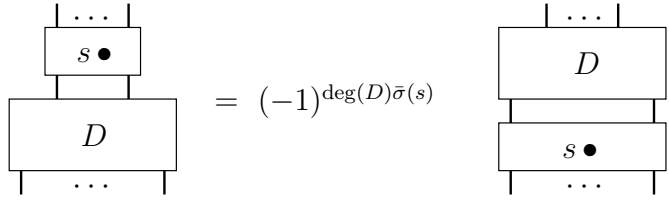}}
$\vcenter{\box1}$
\caption{A generalization of the Relation~($\sigma$-Nat)}\label{f:sDNat}
\end{figure}
\end{lemma}

A special case of Lemma~\ref{sDNat} is shown in
Figure~\ref{f:sNat2}.

\begin{figure}[!h]
\centering
$
\begin{picture}(2,1)(-1,-0)
\thicklines \put(-1,1){\vector(1,0){2}}
\put(-1,-0.5){\vector(1,0){2}} \thinlines
\put(-0.425,-0.5){\line(1,2){0.75}} \put(0.3,-0.5){\circle*{0.15}}
\put(-0.4,1){\circle*{0.15}}
\put(-0.6,0.4){\makebox(0.4,0.5){$s^{-1}$}}
\put(0.1,-1){\makebox(0.4,0.5){$s$}}
\end{picture}
\ = \ \sigma(s) \
\begin{picture}(2,1.1)(-1,-0.1)
\thicklines \put(-1,1){\vector(1,0){2}}
\put(-1,-0.5){\vector(1,0){2}} \thinlines
\put(0.3,-0.5){\line(-1,2){0.75}} \put(-0.4,-0.5){\circle*{0.15}}
\put(0.3,1){\circle*{0.15}}
\put(0.15,0.4){\makebox(0.4,0.5){$s^{-1}$}}
\put(-0.6,-1){\makebox(0.4,0.5){$s$}}
\end{picture}$
\bigskip \caption{A version of
Relation~\mbox{($\sigma$-Nat)}}\label{f:sNat2}
\end{figure}
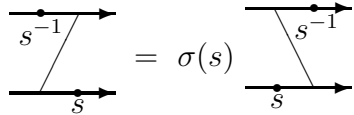

By the Relations~(Rep) we can leave out a label with the neutral
element of a group.
In the case of the projective plane, we represent the non-trivial
element in $\pi_1(P^2\times I,*)\cong\Z/(2)$ by a marking
without label.
We derive the following
particular relations in $\Ab(P^2\times I)$.

\begin{lemma}\label{Z2rels}
In~$\Ab(P^2\times I)$ the relations shown in
Figure~\ref{f:Z2rels} hold, where the markings in the second relation
represent all non-trivial labels of the diagram.

\begin{figure}[!h]
\centering
\setbox1=\hbox{\input{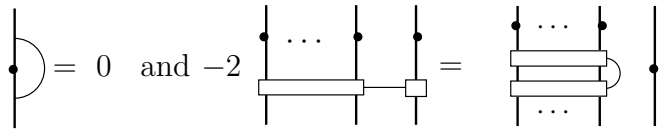}}
$\vcenter{\box1}$
\caption{Special relations for the thickened projective
plane}\label{f:Z2rels}
\end{figure}
\end{lemma}
{\bf Proof:} By using Relation~$(FI)$ and by applying
Lemma~\ref{implrels} in two different ways, we obtain the
equations from Figure~\ref{isochord}.

\begin{figure}[!h]
\centering
\setbox1=\hbox{\input{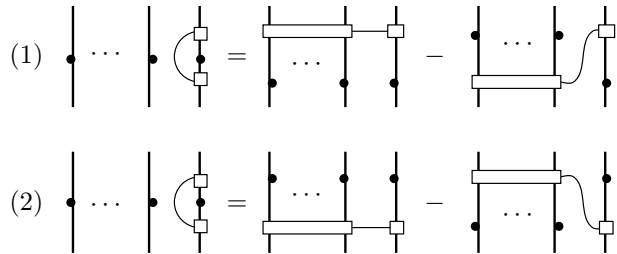}}
$\vcenter{\box1}$
\caption{Proof of a "\labelled{} framing independence" relation}\label{isochord}
\end{figure}

We can apply Lemma~\ref{sDNat} no matter how the strands on the
right side of part~(2) of Figure~\ref{isochord} are oriented
because $x^{-1}=x$ for $x\in\pi_1(P^2\times I,*)$. So when we add
the two equations in Figure~\ref{isochord} the right sides of the
equations cancel and we obtain the first relation of
Figure~\ref{f:Z2rels}. The second relation follows from
Figure~\ref{f:msquare}.

\begin{figure}[!h]
\centering
\setbox1=\hbox{\input{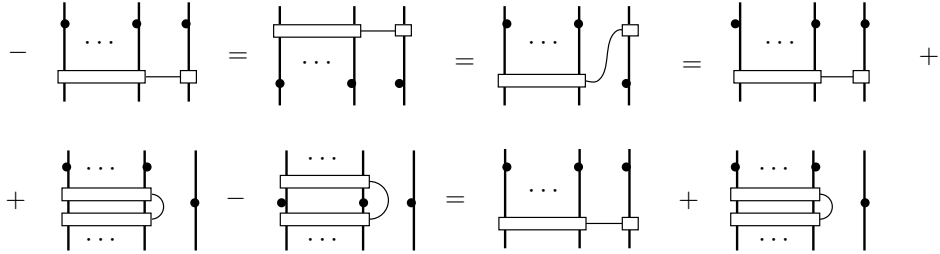}}
$\vcenter{\box1}$
\caption{Freeing the strand on the right side from its
chords}\label{f:msquare}
\end{figure}

In this figure the first equality follows from Lemma~\ref{sDNat}
and the second one follows from the first part of this proof, the
third one follows by~$m-1$ applications of Lemma~\ref{implrels},
and the last one follows by Relation~($\sigma$-Nat), the relation
in Figure~\ref{f:sNat2} and the first part of this proof. $\Box$

 \section{Tangles in a thickened decomposed surface}\label{s:tangles}

Recall that the definition of the map $\psi_M$ only depends on the
choice of a basepoint~$*$ and an oriented neighborhood $*\in
U\cong\R^3\subset M$. The definition of the universal Vassiliev
invariant $Z_{\Sigma\times I}$ will depend on more choices. We
collect the necessary data in the following definition.

\begin{defi}\label{d:decsur}
A {\em decomposed surface} is a surface~$\Sigma$ together with a
tuple of distinguished subsets~$(B_0,B_1,\ldots,B_k,J_1, J_2)$,
where~$I^2\cong B_0\subset\Sigma$ is oriented, $I\cong J_i \subset
B_0\cap\partial\Sigma$ ($i=1,2$), and $I^2\cong B_i\subset\Sigma$
($i>0$) such that $\bigcup_{i=0}^k B_i=\Sigma$ and
equations~(\ref{e:decsur1}) to~(\ref{e:decsur3}) are satisfied.

\begin{eqnarray}
& & B_0\cap B_i\cong I\times\{0,1\}\mbox{ for }i\geq
1,\label{e:decsur1}\\ & & B_i\cap B_j=\emptyset\mbox{ for
}i\not=j, i,j\geq 1 \mbox{ and }J_1\cap
J_2=\emptyset,\label{e:decsur2}\\
 & & \exists J_3\cong I, J_3\subset
B_0\cap\partial\Sigma\mbox{ such that } J_3\cap J_i=\{P_i\}\
(i=1,2)\ \mbox{and}\nonumber\\ & & \mbox{$J_3$ with the
orientation induced by~$\partial B_0$ directs from~$P_1$
to~$P_2$.} \label{e:decsur3}
\end{eqnarray}
\end{defi}

For every connected compact surface~$\Sigma$ there exist subsets
$(B_0,\ldots B_k, J_1, J_2)$ that equip~$\Sigma$ with the
structure of a decomposed surface. The
index~$k=\rank(H_1(\Sigma))$ is uniquely determined by~$\Sigma$.
We call~$B_0$ the disk of~$\Sigma$ and~$B_i$ the bands of the
decomposition. We equip~$J_1$ with the orientation induced
by~$\partial B_0$ and~$J_2$ with the opposite of the orientation
induced by~$\partial B_0$. An isomorphism of decomposed
surfaces~$\Sigma$ and~$\Sigma'$ is a homeomorphism of the surfaces
that respects the additional data.

We represent decomposed surfaces graphically as shown in
Figure~\ref{f:decsur} by an example. By convention the orientation
of~$B_0$ is counterclockwise in this figure. The intervals~$J_i$
direct from left to right. If~$B_0\cup B_i$ is non-orientable,
then the $2$-dimensional representation of the band~$B_i$ is drawn
with a singular point. We will also represent~$\Sigma\times I$ as
in Figure~\ref{f:decsur} where we assume that on~$B_0$ the
interval~$I$ directs towards the reader.

\begin{figure}[!h]
\centering \setbox1=\hbox{\input{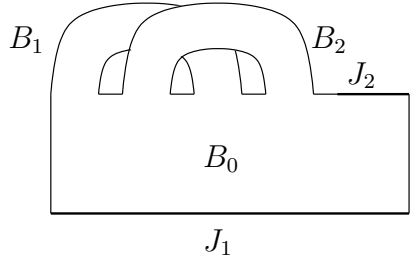}}
$\vcenter{\box1}$ \caption{A picture of a decomposed
surface}\label{f:decsur}
\end{figure}

Given two decomposed surfaces $\Sigma$ and $\Sigma'$
with distinguished subspaces

\begin{equation}
(B_0,\ldots,B_k, J_1,J_2)\ \mbox{and}\
(B'_0,\ldots,B'_{k'}, J'_1,J'_2)
\end{equation}

respectively, we choose an orientation preserving homeomorphism
$\varphi:J_1\longrightarrow J_2'$. We call the surface
$\Sigma\cup_\varphi \Sigma'$ with its natural decomposition

\begin{equation}\label{e:proddecomp}
(B_0\cup B'_0,B_1,\ldots,B_k,B'_1,\ldots,B'_{k'}, J_1',J_2)
\end{equation}

the sum of $\Sigma$ and $\Sigma'$. Let ${\cal S}$ be the set of
isomorphism classes of decomposed surfaces. With the sum of
decomposed surfaces as binary operation ${\cal S}$ becomes a
monoid with neutral element represented by an arbitrary
decomposition of $I^2$. See Figure~\ref{f:bcs} for an example.

\begin{figure}[!h]
\centering \setbox1=\hbox{\input{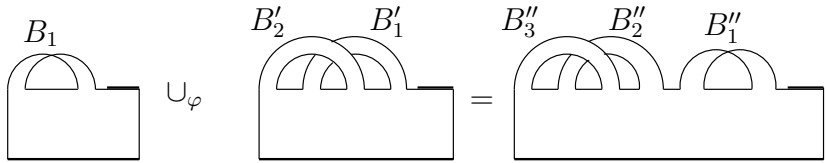}}
$\vcenter{\box1}$ \caption{The sum of decomposed
surfaces}\label{f:bcs}
\end{figure}

\begin{defi}
Let $\Sigma$ be a decomposed surface.
A one-dimensional compact
oriented submanifold~$T\subset\Sigma\times I$ is called a {\em tangle} if

$$
\partial T=T\cap\partial\Sigma\times I\subset \left((J_1\setminus\partial J_1)\cup
(J_2\setminus\partial J_2) \right)\times 1/2.
$$

\end{defi}

We call the ordered sets $T\cap (J_1\times 1/2)$ the lower
boundary and $(T\cap J_2\times 1/2)$ the upper boundary of~$T$. We
assign a symbol~$+$ or~$-$ (resp.~$-$ or~$+$) to a lower (resp.
upper) boundary point of the tangle~$T$ according to whether~$T$
directs towards the boundary point or not. This way we assign
words in the letters~$\{+,-\}$ to the lower and to the upper
boundary of~$T$ called~$\sour(T)$ and~$\tar(T)$ respectively.

Let $\Sigma$ and $\Sigma'$ be decomposed surfaces. We say that two
tangles $T\subset\Sigma\times I$ and $T'\subset\Sigma'\times I$
are {\em isotopic}, if there exists an isotopy $f_t:\Sigma\times
I\longrightarrow\Sigma'\times I$ such that

\begin{eqnarray*}
& & f_0=\alpha\times\id_I, \mbox{where
  $\alpha:\Sigma\longrightarrow\Sigma'$ is an isomorphism of
  decomposed surfaces,}\\
& & {f_t}_{\vert\partial(\Sigma\times I)}={\left(\alpha\times
\id_I\right)}_{\vert\partial(\Sigma\times I)}\mbox{ for all $t\in
[0,1]$,}\\ & & f_1(T)=T'.
\end{eqnarray*}

A link in $\Sigma\times I$ is the same as a tangle
$T\subset\Sigma\times I$ with~$\partial T=\emptyset$.

Let $\Sigma$ be a decomposed surface. A tangle
$T\subset\Sigma\times I$ is called a {\em standard} tangle if it
is contained in $\Sigma\times 1/2$, $T$ has no circle components,
every strand in~$T\cap (B_i\times I)$ ($i>0$) connects the two
components of~$(B_i\cap B_0)\times I$, and every strand in $T\cap
(B_0\times I)$ connects $J_1\times I$ with $(\partial B_0\setminus
J_1)\times I$. See Figure~\ref{f:stdtype} for an example.

\begin{figure}[!h]
\centering
\setbox1=\hbox{\input{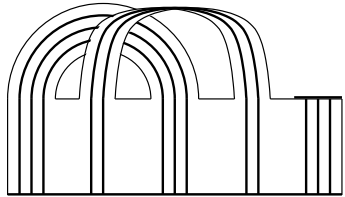}}
$\vcenter{\box1}$
\caption{A standard tangle}\label{f:stdtype}
\end{figure}

Define a category ${\cal T}({\cal S})$ as follows.
The class of objects of ${\cal T}({\cal S})$ is the set of
words~$\{+,-\}^*$
and the morphisms are
isotopy classes of tangles in a decomposed surface.
The composition $T\circ T'$ of tangles $T,T'$ with
$\sour(T)=\tar(T')$ is given by
$T_r\cup T'_r\subset(\Sigma\cup_\varphi\Sigma')\times I$, where $T_r$ and
$T'_r$ are representatives of $T$ and $T'$, such that $T_r\cup T_r'$ is
a tangle in $(\Sigma\cup_\varphi\Sigma')\times I$. The identity
morphisms of ${\cal T}({\cal S})$ are represented by the standard
tangles in~$I^2\times I$.

Let ${\cal T}$ be the subcategory of ${\cal T}({\cal S})$ that has
the same objects as ${\cal T}({\cal S})$, but whose morphisms are
only the isotopy classes of tangles in~$I^2\times I$. We represent
morphisms of ${\cal T}$ always as tangles in $[a,b]\times I\times
I$ for suitable numbers \mbox{$a<b$}, where the surface
$[a,b]\times I$ has the decomposition $([a,b]\times I,[a,b]\times
0,[a,b]\times 1)$. The category ${\cal T}$ is a strict tensor
category. The tensor product of objects is given by the
concatenation of words. The union of $T_1\subset [a,b]\times
I\times I$ with $T_2\subset [b,c]\times I\times I$ induces the
tensor product $T_1\btimes T_2$ of morphisms $T_i$, where we
denote the tensor product of~${\cal T}$ by~$\btimes$ to avoid
confusion with the tensor product~$\otimes$ of vector spaces.
 \section{The functor $Z$}\label{s:recallZ}

In this section we recall the definition and some properties
of a functor~$Z$ from tangles in~$I^2\times I$ to chord diagrams.
Some properties of this functor can be formulated most naturally by
using so-called non-associative tangles that we introduce first.

Given a set $X$, the set of non-associative words $X^{*\na}$ in $X$
is defined as the smallest set such that
the empty word~$\emptyset$ is a non-associative word,
each element~$a\in X$ is a non-associative word,
and if~$r_1$ and~$r_2$ are non-empty non-associative words, then
so is~$(r_1 r_2)$.
We also use the notations $(\emptyset r)$ and $(r\emptyset)$ for the
non-associative word~$r$.
For a non-associative word~$r$ we define the  word~$\bar{r}\in X^*$
by forgetting the brackets of~$r$.
For a category~${\cal C}$ whose objects are words in~$X$, we denote
by~${\cal C}^\na$ the category whose objects are non-associative
words in~$X$ and whose sets of morphisms are given by

\begin{equation}
\Hom_{{\cal C}^\na}(r_1,r_2)=\Hom_{\cal C}(\bar{r_1},\bar{r_2}).
\end{equation}

If ${\cal C}$ is a strict tensor category whose tensor product is
given for objects by the concatenation of words, then one defines
a functor $\btimes:{\cal C}^\na\times{\cal C}^\na\longrightarrow
{\cal C}^\na$ by defining $r_1\btimes r_2=(r_1r_2)$ for objects
of~${\cal C}^\na$ and by defining the tensor product of morphisms
as in~${\cal C}$. We can define associativity constraints
$a_{r_1,r_2,r_3}\in\Hom_{{\cal C}^\na}((r_1r_2)r_3),(r_1(r_2r_3))$
by $a_{r_1,r_2,r_3}=\id_{\bar{r_1}\bar{r_2}\bar{r_3}}$. Then
$({\cal C}^\na,\btimes,\emptyset,a,\id,\id)$ is a tensor category
(see Definition~XI.2.1 of~\cite{Kas}). Starting from the
category~${\cal T}$ of tangles
we consider the category~${\cal T}^\na$ of {\em non-associative tangles}
as a tensor category in this way.

Let $\widehat{\A}({\cal S})$ be the category whose objects are elements
of~$\{+,-\}^*$ and whose morphisms in $\Hom_{\widehat{\A}({\cal S})}(r_1,r_2)$
are (k+1)-tuples $(y,\sigma_1,\ldots,\sigma_k)$ for some $k\in\N$,
where $y\in\widehat{\A}(F_k,\sigma,r_1,r_2)$, $F_k$
is the free group in~$k$ generators $x_1,\ldots,x_k$, and
$\sigma:F_k\longrightarrow\{\pm 1\}$ is given by~$\sigma(x_i)=\sigma_i$.
The composition of two morphisms is given by

$$
(y,\sigma_1,\ldots,\sigma_k)\circ
(z,\sigma'_1,\ldots,\sigma'_{k'})=
(i_*(y)\circ j_*(z),\sigma_1,\ldots,\sigma_k,\sigma'_1,\ldots,\sigma'_{k'}),
$$

where $i:F_k\longrightarrow F_{k+k'}$
and $j:F_{k'}\longrightarrow F_{k+k'}$ are the inclusions given by
$i(x_i)=x_i$ and $j(x_i)=x_{k+i}$.

Let $\widehat{\A}$ be the subcategory of $\widehat{\A}({\cal S})$
that has the same objects as $\widehat{\A}({\cal S})$, but whose
morphisms are only the $1$-tuples $(y)\in\Hom_{\widehat{\A}({\cal
S})}(r_1,r_2)$. The category $\widehat{\A}$ is a strict tensor
category. The tensor product~$\btimes$ of objects is given by the
concatenation of words and the tensor product of morphisms is
induced by the concatenation of their sources and targets.

The category $\widehat{\A}^\na$ has many associativity constraints
turning $\widehat{\A}^\na$ into a tensor category. Recall
from~\cite{Dri} that a Drinfeld associator is a power
series~$\Phi(A,B)$ in two non-commuting indeterminates~$A$ and~$B$
and constant term~$1$ satisfying certain properties. One of these
properties ensures that the morphism $a'_{r_1,r_2,r_3}$
($r_i\in\{+,-\}^{*\na}$) defined in Figure~\ref{f:asso} is an
associativity constraint of~$\widehat{\A}^\na$ and that
$(\widehat{\A}^\na,\btimes,\emptyset,a',\id,\id)$ is a tensor
category. In Figure~\ref{f:asso} the product of indeterminates is
replaced by the composition of chord diagrams as defined in
Section~\ref{s:lcd}.

\begin{figure}[!h]
\centering
\setbox1=\hbox{\input{asso}}
$\vcenter{\box1}$
\caption{Substitution of the non-commuting indeterminates~$A$ and~$B$
of~$\Phi$}\label{f:asso}
\end{figure}

Let $c_{r_1,r_2}\in\Hom_{{\cal T}^\na}((r_1r_2),(r_2r_1))$,
$b_\epsilon\in \Hom_{{\cal T}^\na}(\emptyset,(\epsilon\,-\epsilon))$,
$d_\epsilon\in\Hom_{{\cal T}^\na}((-\epsilon\,\epsilon),\emptyset)$ be
the morphisms of ${\cal T}^\na$ shown in Figure~\ref{f:Tzopf}.

\begin{figure}[!h]
\centering
\setbox1=\hbox{\input{tangds}}
$\vcenter{\box1}$
\caption{Morphisms in ${\cal T}^\na$}\label{f:Tzopf}
\end{figure}

The element $\nu\in\End_{\widehat{\A}}(+)$ is defined as shown in
Figure~\ref{f:defnu}.

\begin{figure}[!h]
\centering
\setbox1=\hbox{\input{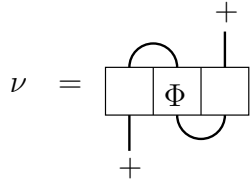}}
$\vcenter{\box1}$
\caption{Definition of the element $\nu$}\label{f:defnu}
\end{figure}

In degree~$0$ the element~$\nu$ is equal to~$\id_+$. Hence there
exists a unique square root~$\nu^{-1/2}$ of~$\nu^{-1}$ that is
equal to~$\id_+$ in degree~$0$. Define~$c'_{r_1,r_2}\in
\Hom_{\widehat{\A}^\na}((r_1r_2),(r_2r_1))$, $b'_\epsilon\in
\Hom_{\widehat{\A}^\na}(\emptyset,(\epsilon\,-\epsilon))$, and
$d'_\epsilon\in\Hom_{\widehat{\A}^\na}((-\epsilon\,\epsilon),\emptyset)$
as shown in Figure~\ref{f:Azopf}, where $\exp(x)=\sum_{n=0}^\infty
x^n/n!$ is regarded as a formal power series.

\begin{figure}[!h]
\centering
\setbox1=\hbox{\input{chordds}}
$\vcenter{\box1}$
\caption{Morphisms in $\widehat{\A}^\na$}\label{f:Azopf}
\end{figure}

Recall the following fact (see Theorem~2 of~\cite{LM2}, Theorem~2 of~\cite{NAT})
that plays a fundamental role in the theory of Vassiliev
invariants.

\begin{fact}\label{f:Z}
There exists a unique functor $Z:{\cal T}^\na\longrightarrow
\widehat{\A}^\na$ preserving the tensor product~$\btimes$ such
that $Z(r)=r$ for objects, and $Z(c_{+,+})=c'_{+,+}$,
$Z(d_\epsilon)=d'_\epsilon$, and
$Z(a_{r_1,r_2,r_3})=a'_{r_1,r_2,r_3}$ for $\epsilon\in\{+,-\}$,
$r,r_i\in\{+,-\}^{*\na}$. We have $Z(b_\epsilon)=b'_\epsilon$ and
$Z(c_{r_1,r_2}^{\pm 1})={c'}_{r_1,r_2}^{\pm 1}$.
\end{fact}

Recall that for an $\{e\}$-\labelled{} chord diagram~$D$ one can
define a formal linear combination~$[T_D]$ of tangles
in~$I^2\times I$ similar to our definition of~$[L_D]$ in
Section~\ref{s:psi}. Extend~$Z$ by linearity to formal linear
combinations, and extend the comultiplication to a map
$\widehat{\Delta}:\Hom_{\widehat{\A}^\na}(r_1,r_2)\longrightarrow
\Hom_{\widehat{\A}^\na}(r_1,r_2)\widehat{\otimes}
\Hom_{\widehat{\A}^\na}(r_1,r_2)$. Using this notation we recall
the following fact (see~\cite{LM3}).

\begin{fact}\label{f:Zprop}
(1) For a chord diagram~$D\in\Hom_{\widehat{\A}^\na}(r_1,r_2)$
with~$\deg(D)=n$ we have

$$
Z([T_D]) =  D+\mbox{terms of degree $>n$.}
$$

(2) For a non-associative tangle~$T\subset I^2\times I$ we have

$$ \widehat{\Delta}(Z(T))=Z(T)\widehat{\otimes} Z(T). $$
\end{fact}

For $r\in\{+,-\}^{*\na}$ define $-r$ by interchanging the symbols
$+$ and $-$, and define~$r^\mid$ by reading the word~$r$ from the
right side to the left side. For a tangle $T\in\Hom_{{\cal
T}^\na}(r_1,r_2)$ we define $T^{-}\in\Hom_{{\cal
T}^\na}(-r_2,-r_1)$ by turning $T$ around a horizontal axis by an
angle~$\pi$, we define $T^\mid\in\Hom_{{\cal
T}^\na}(r_1^\mid,r_2^\mid)$ by turning~$T$ around a vertical axis
by an angle~$\pi$, and we define~$T^\star$ by reflection of~$T$
in~$I^2\times 1/2$. For \mbox{$t=\sum_{i=0}^\infty
t_i\in\Hom_{\widehat{\A}^\na}(r_1, r_2)$} with~$\deg(t_i)=i$ we
define $t^\star\in\Hom_{\widehat{\A}^\na}(r_1, r_2)$ by
\mbox{$t^\star=\sum_{i=0}^\infty(-1)^i t_i$}. We define
$t^{-}\in\Hom_{\widehat{\A}^\na}(-r_2,-r_1)$ and
$t^\mid\in\Hom_{\widehat{\A}^\na}(r_1^\mid,r_2^\mid)$ like~$T^{-}$
and~$T^\mid$ by applying the respective rotations to the graphical
representation of chord diagrams.

A Drinfeld associator~$\Phi$ is called even if it is~$0$ in odd
degrees. By~Proposition~5.4 of \cite{Dri} there exists an even
associator with rational coefficients that is an exponential of a
Lie series. We will fix a choice of an associator with these
properties until the end of Section~\ref{s:s1s2}. With this choice
for the definition of the associativity constraint of the tensor
category~$\widehat{\A}^\na$ the invariant~$Z$ has the following
symmetry properties (see~Proposition~3.1 of~\cite{LM3}).

\begin{fact}\label{f:assosym}
For an associator~$\Phi$ as chosen above one has

$$ Z(T^{-})=Z(T)^{-}\quad \mbox{ , }\quad
Z(T^\mid)=Z(T)^\mid\quad\mbox{and}\quad Z(T^\star)=Z(T)^\star. $$
\end{fact}
 \section{Pairs and diagrams of tangles in $\Sigma\times I$}\label{s:pd}

The existence of the functor~$Z$ of the previous section can be
proved by using a presentation of the category~${\cal T}$ by
generators and relations. In this section we extend this
presentation to~${\cal T}({\cal S})$.
Let~$\Sigma$ be a decomposed
surface. Let~$T$ be a tangle in~$\Sigma\times I$. Then~$T$ is
isotopic to the composition~$T_1\circ T_2$ of a standard
tangle~$T_1\subset \Sigma\times I$ and a tangle~$T_2\subset
I^2\times I$.

\begin{defi}
(1) A pair $(T_1,T_2)$ of a standard tangle $T_1\subset\Sigma\times I$
and of a
tangle $T_2\subset I^2\times I$ with $\sour(T_1)=\tar(T_2)$
is said to be a {\em representing pair}
of the tangle $T_1\circ T_2\subset\Sigma\times I$.

(2) Two representing pairs $(T_1,T_2)$ and $(T_1',T_2')$ are {\em
  isotopic} if $T_1$ is isotopic to $T_1'$ and $T_2$ is isotopic to $T_2'$.
\end{defi}

The graphical representation of the
composition~$T_1\circ T_2$ of a representing pair as
shown in Figure~\ref{f:lidi} by an example is called a {\em diagram}
of a tangle in~$\Sigma\times I$.
Notice that $T_2\subset I^2\times I$ lies in an oriented manifold.
In a diagram, the tangles~$T_1\subset\Sigma\times I$ and
$T_2\subset I^2\times I$ are given by the projection of
representatives in generic position to~$\Sigma$ (resp.\ to $I^2$)
along~$I$ together with the information
saying which one is the overcrossing strand
at double points of the projection of~$T_2$.

\begin{figure}[!h]
\centering
\setbox1=\hbox{\input{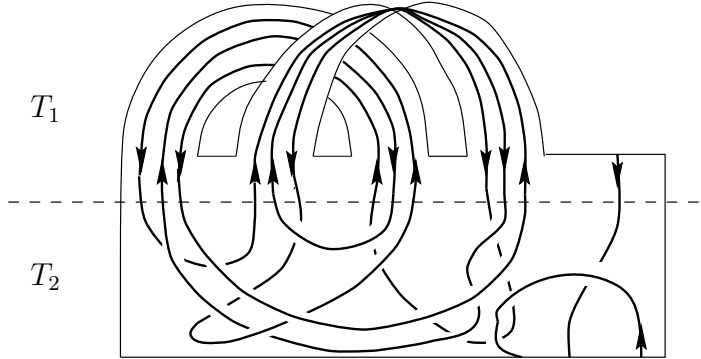}}
$\vcenter{\box1}$
\caption{A diagram of a tangle}
\label{f:lidi}
\end{figure}

In Figure~\ref{f:rmoves} we see four moves between
different representing pairs $(T_1,T_2)$ of a tangle~$T$.
The two parallel strands at the left and at the right side of each
diagram represent a (possibly empty) bundle of strands.
We assume an arbitrary compatible choice
of orientations of the strands.

\begin{figure}[!h]
\centering\setbox1=\hbox{\input{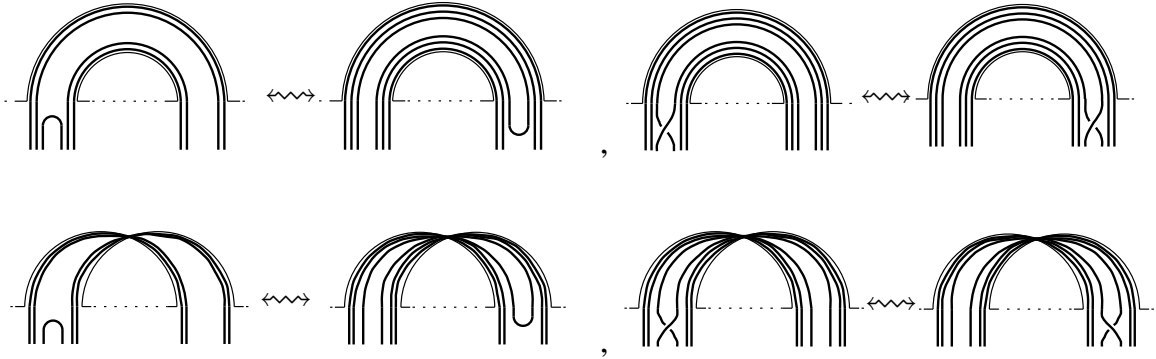}}
$\vcenter{\box1}$ \caption{Moves between representing
pairs of a tangle $T$}\label{f:rmoves}
\end{figure}


Isotopies of tangles in~$\Sigma\times I$ can be described by
sequences of Reidemeister moves and by isotopies of tangle
diagrams.
This can be translated into the following proposition.

\begin{prop}\label{FReidem}
Let $(T_1,T_2)$ and $(T'_1,T'_2)$ be two representing pairs of tangles
in~$\Sigma\times I$. Then $T_1\circ T_2$ is isotopic to
$T'_1\circ T'_2$,
if and only if one can pass from $(T_1,T_2)$ to
$(T_1',T_2')$ by a finite sequence of the following moves:

(1) An isotopy of representing pairs.\nopagebreak

(2) One of the moves shown in Figure~\ref{f:rmoves}, and their images
under reflection in a vertical axis.
\end{prop}
 \section{The functor $Z_{\cal S}$}\label{s:ZS}

In addition to the subcategory ${\cal T}^\na$, the category
${\cal T}({\cal S})^\na$ has a further important
subcategory~${\cal U}$.
The objects of ${\cal U}$ are elements
of~$\{+,-\}^{*\na}$, and its morphisms are standard tangles~$T$ in a
thickened decomposed surface~$\Sigma$
satisfying the following conditions (NA1) and (NA2). To state these
conditions we first introduce some notation:

Let $(B_0,\ldots,B_k,J_1,J_2)$ be the distinguished subsets of the
decomposed surface~$\Sigma$. Let~$T$ be a standard tangle
in~$\Sigma\times I$ with~$s=\sour(T)$. On the upper boundary
of~$B_0$ lie the disjoint intervals $B_i\cap B_0$ and $J_2$. We
denote these intervals by $I_1,\ldots,I_{2k+1}$ numbered from the
left side to the right side. Subsets of the lower boundary of~$T$
correspond to subwords of~$\bar{s}$ and vice versa. The
word~$\bar{s}$ can be subdivided into (possibly empty) subwords
$r_1,\ldots, r_{2k+1}$ such that when we travel along~$T$ starting
at a boundary point of~$T$ belonging to~$r_\nu$ we arrive at the
interval~$I_\nu\times 1/2$ before traveling across any other
interval~$I_\mu\times 1/2$ ($\mu\not=\nu$).

(NA1) There exist non-associative
words $s_\nu$ with $\bar{s_\nu}=r_\nu$ and

$$
s=(s_1(s_2(\ldots(s_{2k-1}(s_{2k}s_{2k+1}))\ldots)))
$$

with all the~$2k$ closing brackets to the extreme right.

(NA2) We have $s_{2k+1}=\tar(T)$. For some~$i\leq k$, let~$I_\nu$
and~$I_\mu$ be the intervals $B_i\cap B_0$. If $B_i\cup B_0$ is
orientable, then we have $s_\nu=-s_\mu^{\mid}$, and $s_\nu=-s_\mu$
otherwise.

\smallskip

For the standard tangle~$T_1$ of the representing pair in
Figure~\ref{f:lidi} an example of parenthesis on~$\sour(T_1)$
satisfying~(NA1) and~(NA2) is the following:

$$ (r_1(r_2(-r_1^{\mid}(-r_2\;+))))=
\left((+(-+))\left(((--)+)\left(((-+)-)\left((++)-)
+\right)\right)\right)\right), $$

where $r_1=(+(-+))$ and $r_2=((--)+)$.

We define a functor $Z_\st:{\cal U} \longrightarrow
\widehat{\A}({\cal S})^\na$ as follows. For objects we have
$Z_\st(r)=r$. For a morphism~$T$ we have
$Z_\st(T)=(z(T),\sigma_1,\ldots,\sigma_k)$, where $\sigma_i=1$ if
$B_0\cup B_i$ is orientable and $\sigma_i=-1$ otherwise, and where
$z(T)$ is the $F_k$-\labelled{} chord diagram of degree~$0$
defined as follows: the skeleton~$\Gamma$ of the chord
diagram~$z(T)$ consists of oriented intervals connecting the same
points as the strands of~$T$. An interval of~$\Gamma$ is
\labelled{} by~$x_i$ (resp.~$x_i^{-1}$) if and only if the
corresponding strand of~$T$ passes through a band~$B_i$ in the
counterclockwise (resp.\ clockwise) sense in a diagram of~$T$.
Using the notation of Figure~\ref{f:sbull} we show an example
for~$z(T)$ in Figure~\ref{f:exaT}.

\begin{figure}[!h]
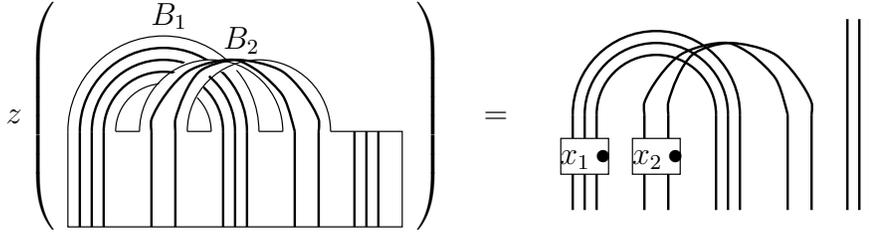

$$
z\left(
\mbox{\setbox1=\hbox{\input{defz}}
$\vcenter{\box1}$}
\right)
\quad =\quad
\mbox{\setbox1=\hbox{\input{defzp2}}
$\vcenter{\box1}$}
$$
\caption{The definition of $z(T)$}\label{f:exaT}
\end{figure}

\begin{theorem}\label{t:TKontsevi}
There exists a unique functor
$Z_{\cal S}:{\cal T}({\cal S})^\na\longrightarrow \widehat{\A}({\cal
  S})^\na$
that makes the following diagram commutative, where the
functors~$i_\nu$ are given by inclusion of subcategories.

$$
\begin{picture}(9.5,5.5)(-4.5,-1)
\put(-4,3){\makebox(2,1){${\cal U}$}}
\put(-4,1){\makebox(2,1){${\cal T}({\cal S})^\na$}}
\put(2,1){\makebox(2,1){$\widehat{\A}({\cal S})^\na$}}
\put(-4,-1){\makebox(2,1){${\cal T}^\na$}}
\put(2,-1){\makebox(2,1){$\widehat{\A}^\na$}}
\put(-1.5,1.5){\vector(1,0){3}}
\put(-1.7,-0.5){\vector(1,0){3.4}}
\put(-3,0){\vector(0,1){1}}
\put(3,0){\vector(0,1){1}}
\put(-3,3){\vector(0,-1){1}}
\put(-2.5,3.5){\vector(3,-1){4.3}}
\put(-1,2.8){\makebox(2,1){$Z_\st$}}
\put(-1,1.5){\makebox(2,0.8){$Z_{\cal S}$}}
\put(-1,-0.5){\makebox(2,0.8){$Z$}}
\put(-4,0){\makebox(1,1){$i_2$}}
\put(3.2,0){\makebox(1,1){$i_3$}}
\put(-4,2){\makebox(1,1){$i_1$}}
\end{picture}
$$
\end{theorem}

{\bf Proof of Theorem~\ref{t:TKontsevi}:}
For objects we necessarily have $Z_{\cal S}(r)=r$.
For morphisms~$T$ the definition of
$Z_{\cal S}(T)$ will involve some choices.
Let~$(T_1,T_2)$ be a representing pair of the tangle~$T$.
Choose~$s\in\{+,-\}^{*\na}$ such that~$T_1$ regarded as a
non-associative tangle with $\tar(T_1)=\tar(T)$ and
$\sour(T_1)=s$ is a morphism of ${\cal U}$.
Regard~$T_2$ as a non-associative tangle with $\sour(T_2)=\sour(T)$
and $\tar(T_2)=s$. Define

\begin{equation}
Z_{\cal S}(T)=Z_\st(T_1)\circ Z(T_2).
\end{equation}

We have to show that $Z_{\cal S}$ is well-defined. Therefore we
verify in part~(a) of this proof that $Z_{\cal S}(T)$ does not
depend on the choice of the non-associative structure of~$s$, and
in part~(b) that $Z_{\cal S}$ is compatible with the moves between
representing pairs of a tangle in Proposition~\ref{FReidem}.

\smallskip

(a) As a special case of Fact~\ref{f:assosym} we have the
equalities

$$
a'^{{\mid}{-}}_{r_1,r_2,r_3}=a'_{-r_1^{\mid},-r_2^{\mid},-r_3^{\mid}}\quad\mbox{and}\quad
a'^{-}_{r_1,r_2,r_3}={a'}_{-r_1,-r_2,-r_3}^{-1} $$

where $r_i\in\{+,-\}^{*\na}$.
This implies the two equations of Figure~\ref{f:Phisym},
where the source of the first morphism is

$$ (((r_1r_2)r_3)((-r_3^{\mid}-r_2^{\mid})-r_1^{\mid})) $$

and the source of the
second morphism is $(((r_1r_2)r_3)(-r_1(-r_2-r_3)))$.

\begin{figure}[!h]
\centering \setbox1=\hbox{\input{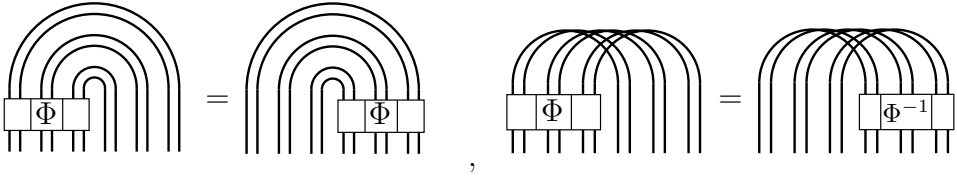}}
$\vcenter{\box1}$ \caption{Symmetries of an even
associator}\label{f:Phisym}
\end{figure}

By Lemma~\ref{sDNat} the equation in Figure~\ref{f:assolabel} holds
because the associator~$\Phi$ is even and therefore the sign in
Figure~\ref{f:sDNat} is always equal to~$1$.

\begin{figure}[!h]
\centering \setbox1=\hbox{\input{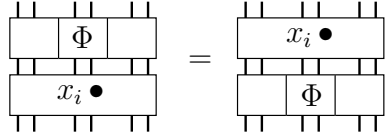}}
$\vcenter{\box1}$ \caption{Commuting an even associator
with the labels $x_i$}\label{f:assolabel}
\end{figure}

Combining the equations from Figure~\ref{f:Phisym}
and~\ref{f:assolabel} we see that for a word~$s$ satisfying~(NA1)
and~(NA2) we can rearrange the parenthesis of~$s_\nu$ and~$s_\mu$
(with the notation of Condition~(NA2)) without changing the value
of $Z_{\cal S}(T)$. But by~(NA1) the parenthesis of the
subwords~$s_\nu$ ($\nu\in\{1,\ldots,2k\}$) of~$s$ are the only
indeterminacy in the choice of~$s$.

\smallskip

(b) The invariance under Move~(1) of Proposition~\ref{FReidem} is
clear because~$Z$ and~$Z_\st$ are well-defined. Let us check the
invariance under the third move of Figure~\ref{f:rmoves}. We can
assume that the parenthesis belonging to the lower boundary of the
shown diagrams have the following properties: Condition~(NA1) is
satisfied, Condition~(NA2) is satisfied for the bands~$B_i$
different from the one shown in the picture, and the source of the
diagrams is of the form $((r_1(\gamma
\;-\gamma))r_2)\ldots(-r_1-r_2)$, where $\gamma\in\{+,-\}$ and
$r_1,r_2\in\{+,-\}^{*\na}$. Then the values of~$Z_{\cal S}$ on the
parts of the diagrams are shown in Figure~\ref{f:im2b}.

\begin{figure}[!h]
\centering
\setbox1=\hbox{\input{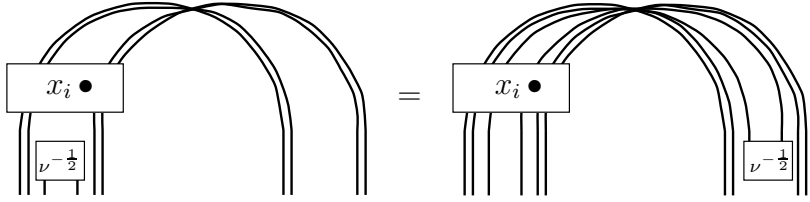}}
$\vcenter{\box1}$
\caption{Invariance under the third move of
Figure~\ref{f:rmoves}}\label{f:im2b}
\end{figure}

By definition~$\nu$ and therefore also $\nu^{-1/2}$
vanish in odd degrees.
Now Lemma~\ref{sDNat} implies that the two diagrams in
Figure~\ref{f:im2b} represent equal elements.

With the assumption on the chosen bracketing similar to the one made
before, the values of~$Z_{\cal S}$ on the parts of the
diagrams in the fourth move of Figure~\ref{f:rmoves}
are shown in Figure~\ref{f:im3b}

\begin{figure}[!h]
\centering \setbox1=\hbox{\input{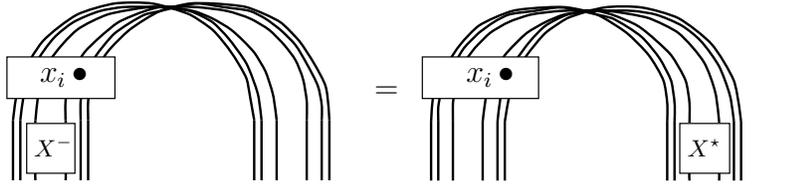}}
$\vcenter{\box1}$ \caption{Invariance under the fourth move
of Figure~\ref{f:rmoves}}\label{f:im3b}
\end{figure}

where $X=\exp\left(\mbox{\setbox1=\hbox{\begin{picture}(0,0)%
\includegraphics{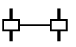}%
\end{picture}%
\setlength{\unitlength}{3947sp}%
\begingroup\makeatletter\ifx\SetFigFont\undefined%
\gdef\SetFigFont#1#2#3#4#5{%
  \reset@font\fontsize{#1}{#2pt}%
  \fontfamily{#3}\fontseries{#4}\fontshape{#5}%
  \selectfont}%
\fi\endgroup%
\begin{picture}(329,200)(1771,-2072)
\end{picture}
}
$\vcenter{\box1}$}/2\right)$. Again an
application of Lemma~\ref{sDNat} shows that the diagrams of
Figure~\ref{f:im3b} represent equal elements. The proof of the
compatibility with the first two moves of Figure~\ref{f:rmoves} is
similar. The proofs are the same when the diagrams are reflected
in a vertical axis. This implies that~$Z_{\cal S}$ is
well-defined. $\Box$

\medskip

Now we come to the first part of the proof of Lemma~\ref{l:Kontsevi}.

\medskip

{\bf Proof of Lemma~\ref{l:Kontsevi} for
  $\partial\Sigma\not=\emptyset$:}
Let $\Sigma$ be a compact connected surface with non-empty
boundary. We choose a decomposition $(B_i,J_i)$ of~$\Sigma$ as in
Definition~\ref{d:decsur} and~$*\in B_0\times I$.
Let~$g_\nu\in\pi_1(\Sigma\times I,*)$ be the element given by a
path that passes one time through the band~$B_\nu$ in the
counterclockwise sense in our pictures. Let
$\varphi:F_k\longrightarrow\pi_1(\Sigma\times I,*)$ be the
isomorphism mapping~$x_\nu$ to~$g_\nu$. Let~$L$ be a link
in~$\Sigma\times I$. We can regard~$L$ as an element
of~$\End_{{\cal T}({\cal S})^\na}(\emptyset)$. Then we define

$$
Z_{\Sigma\times I}(L)=\varphi_*\left(Z_{\cal S}(L)\right),
$$

where the map~$\bar{\cal F}(\varphi)=\varphi_*$ defined in
Section~\ref{s:lcd} is extended to completions. For the definition
of $[L_D]$ we choose the oriented neighborhood~$U$ of~$*$ inside
of~$B_0\times I$. Then equation~(\ref{e:univ}) follows from
part~(1) of Fact~\ref{f:Zprop}, and equation~(\ref{e:groupl})
follows from part~(2) of Fact~\ref{f:Zprop} and from the
definition of~$\Delta(D)$ for diagrams~$D$ of degree~$0$. $\Box$
 \section{The invariant $Z_{P^2\times I}$}\label{s:rp2}

Let $P^2$ be the real projective plane. We choose~$B,\X\subset
P^2$, where $\X$ is a M\"obius strip, $B\cong I^2$, $B\cup
\X=P^2$, and $B\cap \X=\partial B=\partial \X$. By an isotopy we
can push every link in $P^2\times I$ into $\X\times I$. We choose
a decomposition $(B_0,B_1,J_i)$ of~$\X$. Then the diagrams of
links in~$X\times I$ shown in Figure~\ref{RP2move} represent
isotopic links in~$P^2\times I$.

\begin{figure}[!h]
\centering \setbox1=\hbox{\input{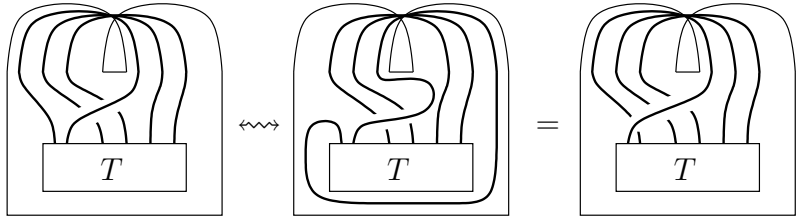}}
$\vcenter{\box1}$ \caption{Pushing a strand of a link
across the disk $B\times 1/2$}\label{RP2move}
\end{figure}

The two strands leaving the box labeled~$T$ in
Figure~\ref{RP2move} represent a possibly empty bunch of strands.
The two diagrams on the right hand side of Figure~\ref{RP2move}
represent isotopic links in~$\X\times I$ (see
Proposition~\ref{FReidem}). A general position argument implies
the following lemma.

\begin{lemma}\label{RP2isotop}
Two links in $\X\times I$ represent isotopic links in $P^2\times I$ if
and only if one can pass from one to the other
by isotopies in $\X\times I$ and by the
move shown in Figure~\ref{RP2move}.
\end{lemma}

The key result for the completion of the proof of Lemma~\ref{l:Kontsevi} is
the following lemma.

\begin{lemma}\label{l:rp2commdiag}
Let $\varphi:\pi_1(\X\times I,*)\longrightarrow\pi_1(P^2\times I,*)$
and~$p:\L(\X\times I)\longrightarrow\L(P^2\times I)$
be the maps induced by the inclusion $\X\times I\subset P^2\times I$.
Then there exists a unique map
$Z_{P^2\times I}:\L(P^2\times I)\longrightarrow
\widehat{\A}(P^2\times I)$
that makes the following diagram commutative.

\nopagebreak

$$
\begin{picture}(9.5,4.5)(-4.5,-2)
\put(-4,1){\makebox(2,1){$\L(\X\times I)$}}
\put(2,1){\makebox(2,1){$\widehat{\A}(\X\times I)$}}
\put(-4,-2){\makebox(2,1){$\L(P^2\times I)$}}
\put(2,-2){\makebox(2,1){$\widehat{\A}(P^2\times I)$}}
\put(-1.2,1.5){\vector(1,0){2.4}}
\put(-1,-1.5){\vector(1,0){2}}
\put(-3,1){\vector(0,-1){2}}
\put(3,1){\vector(0,-1){2}}
\put(-1,1.5){\makebox(2,1){$Z_{\X\times I}$}}
\put(-1,-1.5){\makebox(2,1){$Z_{P^2\times I}$}}
\put(-4,-0.5){\makebox(1,1){$p$}}
\put(3,-0.5){\makebox(1,1){$\varphi_*$}}
\end{picture}
$$
\end{lemma}
{\bf Proof:} We divide each diagram in Figure~\ref{RP2move} into
two parts by cutting along a horizontal line directly above the
box labeled~$T$. Let~$T_1$ (resp.\ $T_2$) be the upper part of the
diagram on the left (resp.\ right) side in Figure~\ref{RP2move}.
We turn~$T_1$ and~$T_2$ into non-associative tangles by
choosing~$\sour(T_1)=\sour(T_2)=(((-\epsilon\, \epsilon)r)-r)$
with appropriate~$\epsilon\in\{+,-\}$ and~$r\in\{+,-\}^{*\na}$.
The values~$t_i=\varphi_*\circ Z_{\X\times I}(T_i)$ are shown in
Figure~\ref{Zvals}.

\begin{figure}[!h]
\centering \setbox1=\hbox{\input{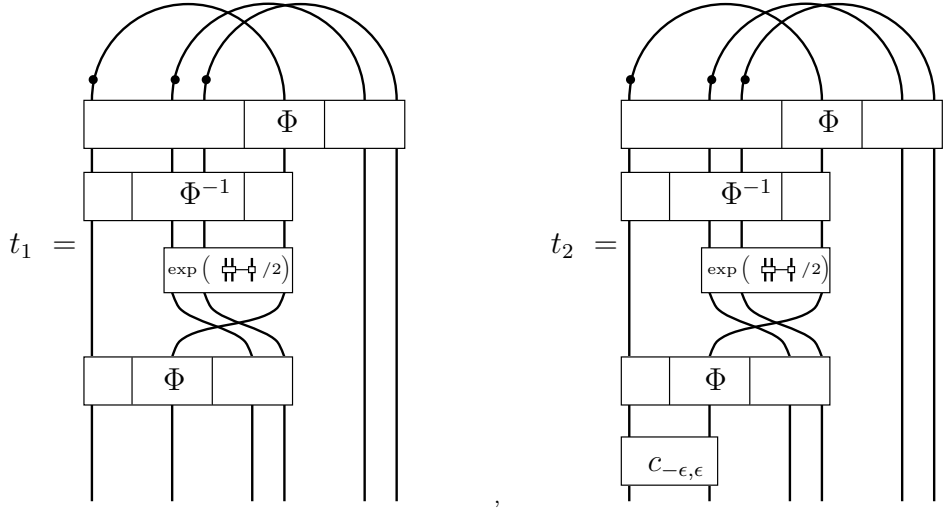}}
$\vcenter{\box1}$ \caption{The values
$t_i=\varphi_*\circ Z_{\X\times I}(T_i)$}\label{Zvals}
\end{figure}

We first consider the associator on the top of these diagrams.
The noncommuting indeterminates~$A$ and~$B$ of~$\Phi$
are replaced by the chord diagrams \setbox1=\hbox{\begin{picture}(0,0)%
\epsfig{file=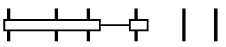}%
\end{picture}%
\setlength{\unitlength}{0.00083300in}%
\begingroup\makeatletter\ifx\SetFigFont\undefined%
\gdef\SetFigFont#1#2#3#4#5{%
  \reset@font\fontsize{#1}{#2pt}%
  \fontfamily{#3}\fontseries{#4}\fontshape{#5}%
  \selectfont}%
\fi\endgroup%
\begin{picture}(1049,200)(1176,-2072)
\end{picture}
}
$\vcenter{\box1}$ and
\setbox1=\hbox{\begin{picture}(0,0)%
\epsfig{file=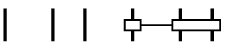}%
\end{picture}%
\setlength{\unitlength}{0.00083300in}%
\begingroup\makeatletter\ifx\SetFigFont\undefined%
\gdef\SetFigFont#1#2#3#4#5{%
  \reset@font\fontsize{#1}{#2pt}%
  \fontfamily{#3}\fontseries{#4}\fontshape{#5}%
  \selectfont}%
\fi\endgroup%
\begin{picture}(1066,200)(1187,-2072)
\end{picture}
}
$\vcenter{\box1}$, respectively.
The equalities in Figure~\ref{firstasso}
follow by Lemma~\ref{implrels} and Figure~\ref{rels}.

\begin{figure}[!h]
\centering
\setbox1=\hbox{\input{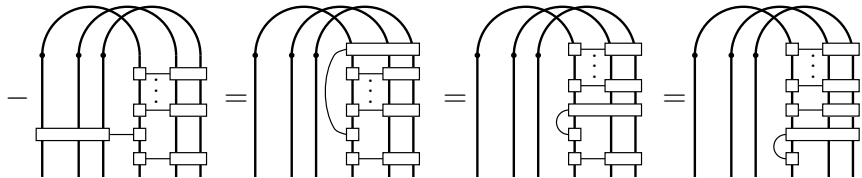}}
$\vcenter{\box1}$
\caption{Replacing the indeterminates in the first
associator}\label{firstasso}
\end{figure}

We see by this figure
that a monomial in~$A$ and~$B$ is replaced by a diagram that only
depends on the number of occurences of~$A$ and~$B$ in the monomial.
But when~$A$ and~$B$ are replaced by commuting elements then the
result vanishes in degree~$>0$ because
the associator~$\Phi$ is an exponential of
a Lie series in~$A$ and~$B$.
We can show by a similar argument using Lemma~\ref{Z2rels}
that the associator in the second box counting from the top in
Figure~\ref{Zvals} vanishes in degree~$>0$ (see
Figure~\ref{f:bundlestrcr} for the replacement of the
indeterminate~$B$ of~$\Phi^{-1}$).

\begin{figure}[!h]
\centering
\setbox1=\hbox{\input{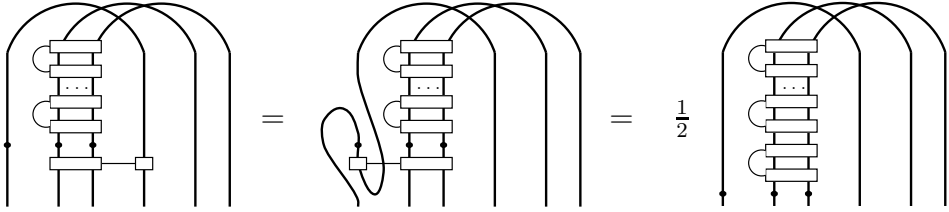}}
$\vcenter{\box1}$
\caption{Freeing the strand on the left side from its
chords}\label{f:bundlestrcr}
\end{figure}

We can calculate the contribution of the third box
of Figure~\ref{Zvals}
using Figure~\ref{f:bundlestrcr}.
It is easy to see that the associator in the fourth
box of Figure~\ref{Zvals} vanishes in degree~$>0$.
Then Lemma~\ref{Z2rels} implies that the contribution from the last
box of~$t_2$ in Figure~\ref{Zvals} also vanishes in degree~$>0$.
Hence we have
the desired equality shown in Figure~\ref{T1gT2}.

\begin{figure}[!h]
\centering
\setbox1=\hbox{\input{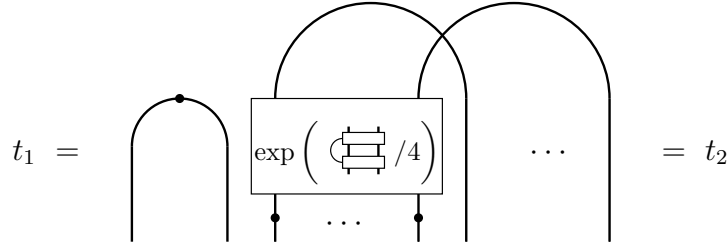}}
$\vcenter{\box1}$
\caption{Invariance of $\varphi_*\circ Z_{\X\times I}$
under the move of Figure~\ref{RP2move}}\label{T1gT2}
\end{figure}

By Lemma~\ref{RP2isotop} this completes the proof.
$\Box$

\medskip

Now it is easy to complete the proof of Lemma~\ref{l:Kontsevi}.

\medskip

{\bf Proof of Lemma~\ref{l:Kontsevi} for $\Sigma=P^2$:} For a
$\pi_1(P^2\times I,*)$-\labelled{} chord diagram $D$ of degree~$n$
we can choose the singular link~$L_D$ lying inside of~$\X\times
I$. The property

$$
Z_{\X\times I}([L_D])=D+\mbox{terms of degree $>n$}
$$

implies the same property for $Z_{P^2\times I}$ in view of
Lemma~\ref{l:rp2commdiag}.
Since $\varphi_*$ is a morphism of coalgebras
the equation

$$
\widehat{\Delta}(Z_{\X\times I}(L))=
\widehat{\Delta}(Z_{\X\times I}(L))\widehat{\otimes}
\widehat{\Delta}(Z_{\X\times I}(L))
$$

implies the same equation for $Z_{P^2\times I}$. $\Box$
 \section{The invariant $Z_{S^1\times S^2}$}\label{s:s1s2}

Let $B\subset S^2$, $B\cong I\times I$. Choose a decomposition
$(B_0,B_1,J_i)$ of the surface $S^1\times I$. We can represent
links in $S^1\times B\cong (S^1\times I)\times I$ by diagrams. The
following lemma follows from general position arguments.

\begin{lemma}\label{kirby}
A link in $S^1\times S^2$ can be represented by a link
in $S^1\times B$. Two links in
$S^1\times B$ represent isotopic links in $S^1\times S^2$ if
and only if they are
the same by isotopies in $S^1\times B$ and by the move shown
in Figure~\ref{KirbyS1S2}.
\end{lemma}

\begin{figure}[!h]
\centering
\setbox1=\hbox{\input{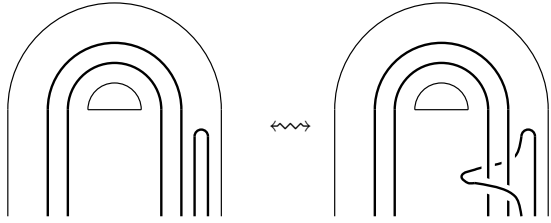}}
$\vcenter{\box1}$
\caption{The second Kirby move for links in $S^1\times S^2$}\label{KirbyS1S2}
\end{figure}

The number of the strands passing through $B_1\times I$ and the
orientations of the strands in Figure~\ref{KirbyS1S2} are
arbitrary. Recall the definition of ${\cal E}(S^1\times S^2)$ from
Section~\ref{s:Z}.

\begin{lemma}\label{l:Kirb2}
The map $\realize_{S^1\times S^2}$ factors through a map

$$
\widetilde{\realize}_{S^1\times S^2}:
{\cal E}(S^1\times S^2)
\longrightarrow\gr\L(S^1\times S^2).
$$
\end{lemma}
{\bf Proof:} Assume that the diagrams in Figure~\ref{KirbyS1S2}
are parts of a singular link having $n-1$ double points and assume
that the strands passing through~$B_1\times I$ in
Figure~\ref{KirbyS1S2} direct from the right to the left. Then the
difference of the desingularizations of the two isotopic singular
links of Figure~\ref{KirbyS1S2} can be expressed as a linear
combination of desingularizations of singular links with~$n$
singularities as shown in Figure~\ref{singlico}.

\begin{figure}[!h]
\centering
\setbox1=\hbox{\input{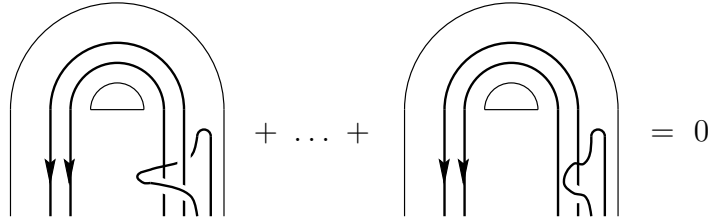}}
$\vcenter{\box1}$
\caption{A relation in $\gr\L(S^1\times S^2)$}\label{singlico}
\end{figure}

Since every part of the link may be pulled to the position of the
maximum on the right hand side of Figure~\ref{KirbyS1S2} it
follows that the Relation~($S^2$-slide) with all strands pointing
upwards is mapped to~$0$ by~$\psi_{S^1\times S^2}$. Considering
what happens when the orientations of some strands in
Figure~\ref{singlico} are changed we see that all
Relations~($S^2$-slide) are mapped to~$0$ by~$\psi_{S^1\times
S^2}$. $\Box$

\medskip

The next lemma will imply that the Relations~($S^2$-slide) generate
{\em all} additional relations in the case $M=S^1\times S^2$.

\begin{lemma}\label{l:ZS1S2}
There exists a unique map $Z_{S^1\times S^2}$ that makes the
following diagram commutative, where~$q:\L(S^1\times I\times I)
\longrightarrow \L(S^1\times S^2)$ and
$\varphi:\pi_1(S^1\times I\times I,*)\longrightarrow \pi_1(S^1\times
S^2,*)$
are induced by the inclusion $S^1\times I\times I\cong S^1\times B
\subset S^1\times S^2$ and $p$ denotes the canonical projection.

$$
\begin{picture}(11.5,5.5)(-5.5,-3)
\put(-5,1){\makebox(2,1){$\L(S^1\times I\times I)$}}
\put(3,1){\makebox(2,1){$\widehat{\A}(S^1\times I\times I)$}}
\put(-5,-3){\makebox(2,1){$\L(S^1\times S^2)$}}
\put(3,-1){\makebox(2,1){$\widehat{\A}(S^1\times S^2)$}}
\put(3,-3){\makebox(2,1){$\widehat{\cal E}(S^1\times S^2)$}}
\put(-2,1.5){\vector(1,0){4}}
\put(-2,-2.5){\vector(1,0){4}}
\put(-4,1){\vector(0,-1){3}}
\put(4,1){\vector(0,-1){1}}
\put(4,-1){\vector(0,-1){1}}
\put(-1,1.5){\makebox(2,0.9){$Z_{(S^1\times I)\times I}$}}
\put(-1,-2.5){\makebox(2,0.9){$Z_{S^1\times S^2}$}}
\put(-5,-1){\makebox(1,1){$q$}}
\put(4,0){\makebox(1,1){$\varphi_*$}}
\put(4,-2){\makebox(1,1){$p$}}
\end{picture}
$$
\end{lemma}
{\bf Proof:} Turn the two tangles~$T_1$ and~$T_2$ shown in
Figure~\ref{KirbyS1S2} into non-associative tangles by choosing
$\sour(T_1)=\sour(T_2)=(r((-r^*\, \epsilon))-\epsilon))$ with
appropriate~$\epsilon\in\{+,-\}$ and $r\in\{+,-\}^{*\na}$. The
values $t_i=\varphi_*\circ Z_{(S^1\times I)\times I}(T_i)$ are
shown in Figure~\ref{bundlew}, where~$s$ is the generator
of~$\pi_1(S^1\times S^2,*)$ that can be represented by a clockwise
oriented circle in our pictures.

\begin{figure}[hpt]
\centering \setbox1=\hbox{\input{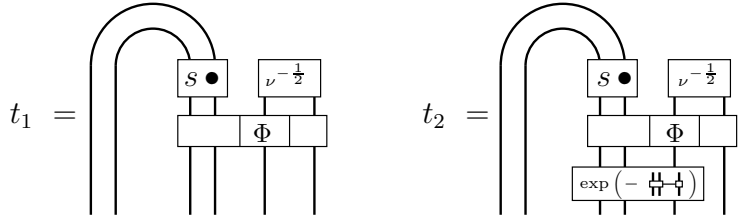}}
$\vcenter{\box1}$ \caption{The values
$t_i=\varphi_*\circ Z_{(S^1\times I)\times I}(T_i)$}
\label{bundlew}
\end{figure}

The associator in the diagrams of $t_1$ and
$t_2$ vanishes in degree $>0$ by the Relation~($S^2$-slide).
Then the expression in the box at the bottom of $t_2$
also vanishes in degree $>0$ by the Relation~($S^2$-slide).
This implies the compatibility with the move shown in
Figure~\ref{KirbyS1S2} and completes the proof.
$\Box$

\medskip

{\bf Proof of Theorem~\ref{t:Kon2}:} The theorem follows by the
same arguments as in the proof of Lemma~\ref{l:Kontsevi} and
Theorem~\ref{th:Kontsevi}. Notice that Proposition~\ref{diagsing}
implies that~$\Ker(\psi_{S^1\times S^2})$ is a coideal. Therefore
${\cal E}(S^1\times S^2)$ is a coalgebra.$\Box$

\medskip

Alternatively to the proof above it may be verified directly that
Relation~($S^2$-slide) generates a coideal of~$\Ab(S^1\times S^2)$
by using Lemma~\ref{implrels}.

Let $S^1\tilde{\times}S^2$ be the quotient of $S^1\times S^2$ by
the fixed point free involution of~$S^1\times S^2$ given
by~$(x,y)\mapsto(-x,-y)$.  Denote by~$p:S^1\times
S^2\longrightarrow S^1\tilde{\times}S^2$ the canonical projection.
Let $\alpha:I\longrightarrow S^1$ and~$\beta:I\longrightarrow S^2$
be paths that connect antipodal points. A closed tubular
neighborhood of~$p((\alpha\times \beta)(I))$ and the closure of
the complement of this neighborhood are both homeomorphic
to~$X\times I$, where~$X$ is the M\"obius strip. Represent links
in $S^1\tilde{\times}S^2$ by links in~$X\times I$. Then the
obvious versions of Lemmas~\ref{l:Kirb2} and~\ref{l:ZS1S2} and of
Theorem~\ref{t:Kon2} hold. This establishes a universal Vassiliev
invariant~$Z_{S^1\tilde{\times}S^2}$ of links
in~$S^1\tilde{\times}S^2$ with values in a space~$\widehat{\cal
E}(S^1\tilde{\times}S^2)$.
 \section{Variations on the definition of $Z_{\Sigma\times
I}$}\label{s:vari}

For some applications the definition of a decomposed
surface~$\Sigma$ is too restrictive. Our first goal in this
section will be to extend the definition of~$Z_{\Sigma\times I}$
to a more general structure on~$\Sigma$. We briefly recall the
notion of a ribbon graph (see~\cite{Tur}, Section~2.1, with the
difference that we allow non-orientable ribbon graphs as well): a
{\em band}~$B$ is homeomorphic to~$I\times I$, has a distinguished
lower base $I\times 0$ and a distinguished upper base $I\times 1$,
and its core $(1/2)\times I$ directs from~$I\times 1$ to~$I\times
0$. A {\em coupon}~$C$ is homeomorphic to~$I\times I$ and has a
distinguished oriented base~$I\times 0$. An (abstract) {\em ribbon
graph} is a decomposition of a surface~$\Sigma$ into bands and
coupons as shown in Figure~\ref{f:ribbonex} by an example.

\begin{figure}[!h] \centering
\setbox1=\hbox{\input{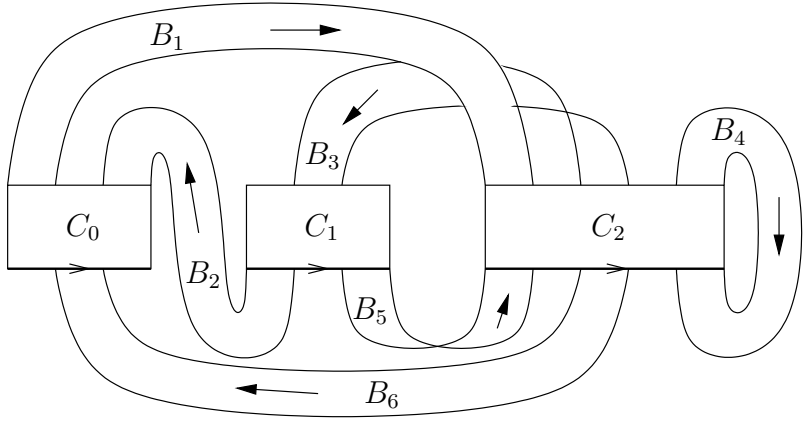}}
$\vcenter{\box1}$ \caption{An example of a ribbon graph}
\label{f:ribbonex} \end{figure}


We consider a locally trivial projection~$p:\Sigma'\times
I\longrightarrow \Sigma$ with typical fiber~$I$, where~$\Sigma'$
is a surface and~$\Sigma$ is a ribbon graph with bands~$B_i$
($i\in J_1$), coupons~$C_j$ ($j\in J_2$), a basepoint~$*\in
C_o\subset\Sigma$ for a distinguished coupon~$C_o$ ($o\in J_2$),
and where~$p^{-1}(C_o)$ is oriented. We choose a maximal number of
bands $B_k$ ($k\in J_3\subset J_1$) such that

\begin{equation}\label{e:chooseT} T=\bigcup_{i\in J_3}
B_i\cup\bigcup_{j\in J_2} C_j \end{equation}

is contractible. We extend the orientation of~$p^{-1}(C_o)$ to
$p^{-1}(T)$. The core of each band~$B$ leads from~$T$ to~$T$ and
therefore defines an element $g_B\in\pi_1(\Sigma,*)$. For $i\in
J_3$ the element $g_{B_i}$ is the neutral element.

Orient~$C_j$ ($j\in J_2$) such that the induced orientation
of~$\partial C_j$ coincides with the orientation of the
distinguished base of~$C_j$. We choose orientation preserving
homeomorphisms~$\varphi_j:C_j\times I\longrightarrow p^{-1}(C_j)$
such that $p(\varphi(c,t))=c$ for all $c\in C_j$, $t\in I$. Let
$L\subset \Sigma'\times I$ be a link. We may assume that
$L\subset\Sigma'\times I$ is in standard position, meaning that
$L\cap p^{-1}(C_j)$ is a tangle in~$p^{-1}(C_j)\cong C_j\times I$
for every $j\in J_2$, $p$ maps $L\cap p^{-1}(B_i)$ injectively
onto $p(L)\cap B_i$ and every component of $p(L)\cap B_i$ connects
the two different bases of the band~$B_i$ ($i\in J_1$).

For a basepoint $*'\in p^{-1}(*)$ the map $p_*:\pi_1(\Sigma'\times
I,*')\longrightarrow \pi_1(\Sigma,*)$ is an isomorphism. To an
oriented strand of $p(L)\cap B_i$ we associate the label
$p_*^{-1}(g_{B_i})$ if its orientation coincides with the
orientation of the core of~$B_i$, and the
label~$p_*^{-1}(g_{B_i}^{-1})$ if the orientations do not
coincide.

Let~$J$ be the distinguished oriented base of a coupon~$C$
of~$\Sigma$. As in Section~\ref{s:tangles} we
associate~$\sour(L\cap p^{-1}(C))\in\{+,-\}^*$ to~$p(L)\cap J$
and~$\tar(L\cap p^{-1}(C))\in\{+,-\}^*$ to~$p(L)\cap \partial
C\setminus J$, where~$\partial C$ carries the opposite of the
orientation induced by~$J$. The conditions~(NA1) and~(NA2) can be
adapted to ribbon graphs in an obvious way, giving us
non-associative tangles~$L\cap p^{-1}(C_j)$. We define
$Z^p_{\Sigma'\times I}(L)$ as in Section~\ref{s:ZS} by gluing
bunches of labeled strands corresponding to~$L\cap p^{-1}(B_i)$ to
the invariants~$Z(L\cap p^{-1}(C_j))$.

\begin{prop}\label{p:Zp} The definition
of~$Z^p_{\Sigma'\times I}$ does not depend on the choice
of~$T\subset\Sigma$. The map~$Z^p_{\Sigma'\times I}$ is an isotopy
invariant of links~$L\subset\Sigma'\times I$. \end{prop} {\bf
Proof:} In order to verify that~$Z_{\Sigma\times I}$ does not
depend on the choice of~$T$, it is enough to consider~$T'$ that
coincides with~$T$ except for two bands that are incident to the
same coupon. This case follows by a single application of
Lemma~\ref{sDNat}.

The rest of the proof proceeds along the same lines as the proof
of Theorem~\ref{t:TKontsevi}. We discuss two differences:

(1) We have to consider bands~$B$ that connect the bottom of a
coupon with the top of a coupon.

The independence of $Z^p_{\Sigma'\times I}(L)$ from the choice of
parentheses coming from the intersection of $p(L)$ with the upper
and lower boundary of~$B$ is simpler because no symmetry with
respect to rotation around a horizontal axis is needed. The
compatibility with moving crossings and local minima or maxima
along~$B$ follows like in the proof of Theorem~\ref{t:TKontsevi}
where we again use less symmetry properties of the elements
$c'_{\epsilon_1\;\epsilon_2}$, $b'_\epsilon$, and~$d'_\epsilon$ in
Figure~\ref{f:Azopf} than before.

(2) Let $\sigma:\pi_1(\Sigma,*)\longrightarrow \{\pm 1\}$ and
$\sigma':\pi_1(\Sigma'\times I,*')\longrightarrow \{\pm 1\}$ be
the orientation characters of~$\Sigma$ and~$\Sigma'\times I$,
respectively. We have to consider bands~$B\subset\Sigma$
with~$-1=\sigma(g_B)\not=\sigma'(p_*^{-1}(g_B))=1$.

Since orientations are only involved in the value of~$Z$ on
crossings the same arguments as before imply that the
map~$Z^p_{\Sigma'\times I}$ does not depend on the choice of
parentheses and is compatible with moving local minima and maxima
of a diagram of~$L$ across the band~$B$. The
map~$Z^p_{\Sigma'\times I}$ is also compatible with moving
crossings across the band~$B$ because~$p^{-1}(T\cup B)$ is
orientable, the preimage~$p^{-1}(C_i)$ of each coupon is equipped
with the orientation induced by~$p^{-1}(T)$,
and~$\widehat{\A}(\Sigma'\times I)$ is defined using~$\sigma'$.
$\Box$

\medskip

As before $Z^p_{\Sigma'\times I}$ is a universal Vassiliev
invariant of links in~$\Sigma'\times I$. The
value~$Z^p_{\Sigma'\times I}(L)$ depends non-trivially on the
ribbon graph structure of~$\Sigma$. There exists an extension of
Proposition~\ref{p:Zp} that allows to consider certain
non-continous projections~$p:\Sigma'\times I\longrightarrow
\Sigma$ as well. As an example this allows to represent links
in~$X\times I$ by diagrams on~$S^1\times I$. I do not know how to
obtain a universal Vassiliev invariant of links in the non-trivial
$I$-bundle over~$P^2$. This invariant would give rise to a
universal Vassiliev invariant of links in~$P^3$.

Besides varying the projection of a link~$L\subset\Sigma'\times I$
there are also different possible normalizations of
$Z^p_{\Sigma'\times I}(L)$ for a fixed diagram of~$L$:
let~$s=s_1\ldots s_n\in\{+,-\}^*$ be a word with~$n$ letters. Then
there exists a map
$d_s:\widehat{\A}(G,\sigma,+,+)\longrightarrow\widehat{\A}(G,\sigma,s,s)$
defined as follows: the skeleton of~$d_s(D)$ consists of a bunch
of~$n$ strands. When $s_i=+$ the labels and the orientation of the
$i$-th strand coincide with the labels and orientation of~$D$.
When~$s_i=-$ the labels are the inverse elements of the labels
of~$D$ and the orientation is opposite. The diagram~$d_s(D)$ is
defined by replacing each chord endpoint by the signed sum of~$n$
ways of lifting that endpoint to the new skeleton, where the signs
are determined by Figure~\ref{multichord}.

After the choice of~$T\subset\Sigma$ as in
equation~(\ref{e:chooseT}) we associate an orientation to a
band~$B$ of~$\Sigma$ as follows: Let~$V$ be a small neighborhood
of the lower base of~$B$ such that the orientation of~$p^{-1}(T)$
extends uniquely to~$p^{-1}(T\cup V)$. Let~$C$ be the coupon
of~$\Sigma$ with~$C\cap V\not=\emptyset$. Orient~$C$ such that the
orientation of~$\partial C$ coincides with the orientation of the
base of~$C$. Let~$\varphi:(V\cup C)\times I\longrightarrow
p^{-1}(V\cup C)$ be a homeomorphism such that~$\varphi_{\vert
C\times I}$ is orientation preserving and~$p(\varphi(x,t))=x$ for
all~$x\in V\cup C$, $t\in I$. Then there exists a unique
orientation of~$B\supset V$ such that~$\varphi_{\vert V\times I}$
is orientation preserving. We denote the band~$B$ with the
orientation from above by~${\rm or}_T(B)$.

Let~${\cal B}^{\rm or}(\Sigma)$ be the set of oriented bands
of~$\Sigma$. For $B\in {\cal B}^{\rm or}(\Sigma)$ we denote
by~$B^\star$ the band~$B$ with the opposite orientation. Let~$a$
be a map from~${\cal B}^{\rm or}(\Sigma)$
to~$\widehat{\A}(\{e\},1,+,+)$ satisfying $a(B^\star)=a(B)^\star$.
We denote~$a({\rm or}_T(B))$ simply by~$a_T(B)$. For a
link~$L\subset \Sigma'\times I$ in standard position,
let~$s(B,L)\in\{+,-\}^*$ be the sequence of letters determined by
the intersection of~$p(L)$ with the lower base of~$B$. Then we
define $Z^{p,a}_{\Sigma'\times I}(L)$ in the same way as
$Z^p_{\Sigma'\times I}(L)$ with the only difference that instead
of gluing bunches of labeled intervals corresponding to $L\cap
p^{-1}(B)$ we now glue the product of bunches of labeled intervals
with~$d_{s(B,L)}(a_T(B))$ to the invariants of the non-associative
tangles~$L\cap p^{-1}(C_j)$.

\begin{prop}\label{p:Zpa}
The definition of~$Z^{p,a}_{\Sigma'\times I}$ does not depend on
the choice of~$T\subset\Sigma$. The map~$Z^{p,a}_{\Sigma'\times
I}$ is an isotopy invariant of links~$L\subset\Sigma'\times I$. If
for all~$B\in {\cal B}^{\rm or}(\Sigma_1)$ the elements~$a(B)$
satisfy~$\widehat{\Delta}(a(B))=a(B)\widehat{\otimes}a(B)$, then
we also have

\begin{equation}\label{e:gl}\widehat{\Delta}\left(Z^{p,a}_{\Sigma'\times
I}(L)\right)=Z^{p,a}_{\Sigma'\times
I}(L)\widehat{\otimes}Z^{p,a}_{\Sigma'\times I}(L).\end{equation}
\end{prop} {\bf Proof:} The first two statements of the
proposition follow along the same lines as Proposition~\ref{p:Zp}
and Theorem~\ref{t:TKontsevi} by using the relations in
Figure~\ref{rels}.

\smallskip

By part 2) of Fact~\ref{f:Zprop} we have $\widehat{\Delta}(t_j)
=t_j \widehat{\otimes}t_j$ for $t_j=Z^{p,a}_{\Sigma'\times
I}\left(L\cap p^{-1}(C_j)\right)$ with $j\in J_2$.
Equation~(\ref{e:gl}) follows by applying properties 1)~to 3)
below to the building blocks~$a_T(B_i)$ ($i\in J_1$) and~$t_j$
($j\in J_2$) of~$Z^{p,a}_{\Sigma'\times I}(L)$.

1)
Let~$\kappa:\widehat{\A}(G,\sigma,s,t)\longrightarrow\widehat{\A}(G,\sigma,s',t')$
be a map induced by gluing some~$G$-labeled intervals to boundary
points of a $G$-labeled diagram~$D$ with~$\sour(D)=s$,
$\tar(D)=t$.
Then~$\widehat{\Delta}\left(\kappa(b)\right)=(\kappa\widehat{\otimes}\kappa)\left(\widehat{\Delta}(b)\right)$.

2) Let $b_i\in\widehat{\A}(G,\sigma,s_i,t_i)$ ($i=1,2$)
with~$\widehat{\Delta}(b_i)=b_i\widehat{\otimes} b_i$. Then
$$\widehat{\Delta}(b_1\btimes b_2)=(b_1\btimes
b_2)\widehat{\otimes}(b_1\btimes b_2).$$

3) For $b\in \widehat{\A}(\{e\},1,+,+)$
with~$\widehat{\Delta}(b)=b\widehat{\otimes}b$ and~$s\in\{+,-\}^*$
we have
$$\widehat{\Delta}\left(d_s(b)\right)=d_s(b)\widehat{\otimes}d_s(b).$$

The properties~1) and 2) follow directly from the definitions. We
give a sketch of a proof of property 3): the equation
$\widehat{\Delta}(b)=b\widehat{\otimes}b$ implies that $b=\exp(P)$
for some~$P\in \widehat{\A}(\{e\},1,+,+)$
with~$\widehat{\Delta}(P)=\id_+\widehat{\otimes }
P+P\widehat{\otimes }\id_+$. By Theorem~8 of~\cite{BN1} the
element~$P$ can be written in terms of connected trivalent
diagrams. The descriptions of the extensions of the
maps~$\widehat{\Delta}$ and~$d_s$ to trivalent diagrams imply that
$\widehat{\Delta}\left(d_s(P)\right)=\id_s\widehat{\otimes
}d_s(P)+d_s(P)\widehat{\otimes }\id_s$. This implies property~3)
and completes the proof. $\Box$

\medskip

Invariants of links in~$\Sigma'\times I$ might be helpful in a
definition of invariants of links in~$I$-bundles over closed
surfaces or of invariants of links in closed $3$-manifolds defined
in terms of Heegard splittings. For these purposes it might be
useful to have different normalizations of~$Z^p_{\Sigma'\times I}$
like provided by Proposition~\ref{p:Zpa}.

Besides different possible choices of diagrams and normalizations,
one can also consider different objects inside of~$\Sigma'\times
I$. Extensions of the definition of~$Z^{p,a}_{\Sigma'\times I}$ to
framed (and half-framed) links, tangles and trivalent graphs with
and without orientation can be adapted without problems
from~$\R^3$ (see~\cite{MuO}, \cite{BeS}) to~$\Sigma'\times I$. In
view of an application in~\cite{Li2} we give some details
concerning an adaptation of $Z^{p,a}_{\Sigma'\times I}$ to links
with basepoints. By definition of a link~$L$ with basepoints there
has to be one distinguished point~$P$ on each component of~$L$.
These basepoints have to lie inside of~$p^{-1}(T)$ and must not
leave~$p^{-1}(T)$ during isotopies. On labeled chord diagrams with
basepoints the basepoints have to be disjoint from labels and
chord endpoints. We represent basepoints graphically by the
symbol~$\times$. We define a $\Q$-vector space~$\Ab^b(G,\sigma)$
generated by $G$-labeled chord diagrams with basepoints modulo the
relations (4T), (FI), (Rep), ($\sigma$-Nat) and modulo the
Relation~(Bas1) shown in Figure~\ref{f:Bas1}.

\begin{figure}[!h] \centering
\setbox1=\hbox{\input{bas1}}
$\vcenter{\box1}$ \caption{The Relation (Bas1)} \label{f:Bas1}
\end{figure}

For a $3$-manifold~$M$ we denote the completion
of~$\Ab^b(\pi_1(M,*),\sigma)$ by~$\widehat{\A}^b(M)$. For a
generic diagram of a link with basepoints~$L\subset\Sigma'\times
I$ we define $Z^{p,a}_{\Sigma'\times
I}(L)\in\widehat{\A}^b(\Sigma'\times I)$ in the same way as for
links without basepoints except that we map basepoints of links to
basepoints of chord diagrams. We have the following lemma.

 \begin{lemma} For a link~$L$ with basepoints the element
$Z^{p,a}_{\Sigma'\times I}(L)\in\widehat{\A}^b(\Sigma'\times I)$
is invariant under isotopies of links with basepoints.
\end{lemma} {\bf Proof:} Relation~(Bas1) assures that~$Z^{p,a}_{\Sigma'\times I}(L)$
is compatible with moving basepoints across local minima, local
maxima, crossings and bands~$B\subset T$ of~$\Sigma$ in a generic
diagram of~$L$. The rest of the proof follows by the same
arguments as in the proofs of Propositions~\ref{p:Zp}
and~\ref{p:Zpa}.$\Box$

\medskip

Links and labeled chord diagrams with basepoints are useful to
avoid problems with signs in the definition of invariants of links
in non-orientable $3$-manifolds. The adaptions of~$Z_{P^2\times
I}$ and~$Z_{S^1\tilde{\times} S^2}$ to links with basepoints are
straightforward.
 \section{Universal Vassiliev invariants and coverings}\label{s:cover}

Let~$H, G$ be groups and let $\sigma_H:H\longrightarrow \{\pm 1\}$
and $\sigma_G:G\longrightarrow \{\pm 1\}$ be homomorphisms. Let
$i:H\longrightarrow G$ be a monomorphism satisfying
$\sigma_G(i(h))=\sigma_H(h)$ for all~$h\in H$. We identify~$H$
with the subgroup~$i(H)$ of~$G$. The group~$G$ acts from the right
on the set of left cosets~$H\setminus G$. Assume that the
index~$d=[G:H]$ is finite. Then~$i$ induces a morphism of
coalgebras

\begin{equation}\label{e:istar}
i^*:\Ab(G,\sigma_G)\longrightarrow \Ab(H,\sigma_H)
\end{equation}

that we describe below. Let~$D$ be a $G$-labeled chord diagram.
Replace the skeleton~$\Gamma$ of~$D$ by a $d$-fold
covering~$i^*(\Gamma)$ of~$\Gamma$ that is constructed as follows.
Replace intervals on~$\Gamma$ without labels by a bunch of~$d$
intervals that are in bijection with~$H\setminus G$. These bunches
of intervals are glued according to the permutations
on~$H\setminus G$ induced by the labels~$g\in G$ of~$\Gamma$. The
labeled points of~$\Gamma$ are covered by points on~$i^*(\Gamma)$
in front of and close to these permutations. In order to describe
how to replace the labels and chords of~$D$, we choose a set
${\cal R}=\{g_a\in a\;\vert\; a \in H\setminus G\}$ of
representatives of left cosets. Then we lift each label~$g\in G$
of~$\Gamma$ to~$d$ labels on~$i^*(\Gamma)$, where the label on the
interval belonging to~$a\in H\setminus G$ is~$g_agg_{a\cdot
g}^{-1}\in H$ (or, more precisely, $j(g_agg_{a\cdot g}^{-1})\in
H$, where $j:i(H)\longrightarrow H$ is the inverse of~$i$). We
complete the description of~$i^*(D)$ by lifting the chords
to~$i^*(\Gamma)$ as follows. We replace each chord of~$D$ by a
linear combination of~$d$ terms
by summing over the lifts of the chord endpoints to
intervals belonging to the same coset~$a$ of~$H\setminus G$ with
coeffiecients~$\sigma(g_a)$.

\begin{lemma} The map $i^*$ is well-defined and does not depend on the choice of~${\cal
R}$.
\end{lemma} {\bf Proof:} First we have to verify that~$i^*$ is
compatible with the relations~($\sigma$-Nat), (Rep), (4T),
and~(FI): consider a relation~($\sigma$-Nat) of the form
$D=\sigma(s) D'$ like in Definition~\ref{d:A}. Expand $i^*(D)$
(resp.\ $i^*(D')$) as a sum of~$d$ terms~$t_a$ (resp.\ $t'_a$)
with $a\in H\setminus G$ according to which pair of strands the
chord in the relation is lifted. Then by relation~($\sigma$-Nat)
in~$\Ab(H,\sigma_H)$ and by the signs belonging to the lifted
chords we have $\sigma (g_{a\cdot s}) i^*(t_{a\cdot s})=\sigma(g_a
s g_{a\cdot s}^{-1})\sigma(g_{a})i^*(t'_{a})$ which
implies~$i^*(t_{a\cdot s} )=\sigma(s)i^*(t'_{a})$. We
obtain~$i^*(D)=\sigma(s)i^*(D')$ by summing over~$a\in H\setminus
G$.

For $g,h\in G$ and $a\in H\setminus G$ we have $g_a gh g_{a \cdot
gh}^{-1}=(g_a g g_{a\cdot g}^{-1})(g_{a\cdot g} h g_{a \cdot
gh}^{-1})$. Since~$G$ acts on~$H\setminus G$ we obtain
compatibility of $i^*$ with the relations~(Rep). The argument
that~$i^*$ is compatible with the relations~(4T) and~(FI) is the
same as in the proof that the comultiplication~$\Delta$ (see
equation~(\ref{e:comult})) is compatible with these relations and
is well-known from the case of diagrams without labels.

Define a map~${i^*}'$ in the same way as~$i^*$, but by replacing
one coset representative~$g_a\in{\cal R}$ by~$g'_a=hg_a$
with~$h\in H$. It remains to show that~$i^*={i^*}'$. Let~$D$ be
a~$G$-labeled chord diagram. We compare~$i^*(D)$ with ${i^*}'(D)$
term by term: we can eliminate all new labels~$h, h^{-1}$
in~${i^*}'(D)$ by using relations~(Rep), then ($\sigma$-Nat) for
each chord~$c$ in~${i^*}'(D)$ between strands corresponding
to~$a$, and again~(Rep). This way we obtain the product
of~$\sigma(h)$ from relation~($\sigma$-Nat) with~$\sigma(g'_a)$
from the definition of~${i^*}'$ for all chords~$c$ as above. This
product equals~$\sigma(g_a)$ from the definition of~$i^{*}(D)$.
 $\Box$

\medskip

Let $p:E\longrightarrow B$ be a finite covering of connected
$3$-manifolds. Taking the preimage of links in~$B$ induces a map
$p^*:\L(B)\longrightarrow \L(E)$.

\begin{lemma}\label{l:pstarfilt} The map $p^*$ satisfies
$p^*(\L_n(B))\subseteq\L_n(E)$ for all~$n\geq 0$.  \end{lemma}
{\bf Proof:} When~$L, L'\subset B$ differ by a crossing change,
then~$p^{-1}(L)$ can be obtained from~$p^{-1}(L')$ by changing~$d$
crossings. Therefore the preimage of the desingularization (see
Figure~\ref{f:desing}) of a singular link in~$B$ with~$n$ double
points can be written as a sum of desingularizations of~$n^d$
singular links in~$E$ with suitable local orientations at the
double points. $\Box$

\medskip

By Lemma~\ref{l:pstarfilt} the map~$p^*$ induces a
map~$\gr(p^*):\gr\L(B)\longrightarrow\gr\L(E)$. Choose a basepoint
$\tilde{*}\in p^{-1}(*)$. From now on we will consider the case
$G=\pi_1(B, *)$ and $H=\pi_1(E,\tilde{*})$, where~$\sigma_G$
(resp.\ $\sigma_H$) is the orientation character of~$B$ (resp.\
$E$), and the inclusion~$H\longrightarrow G$ is given
by~$p_*=\pi_1(p)$. We identify $H$ with $p_*(H)\subset G$. Choose
oriented neighborhoods~$U\subset B$ of~$*$ and
$\widetilde{U}\subset p^{-1}(U)$ of~$\tilde{*}$ with~$U,
\widetilde{U}\cong \R^3$ such that $p_{\vert \widetilde{U}}
:\widetilde{U}\longrightarrow U$ is orientation preserving. Then
by the following proposition the maps~$\pi_1(p)^*$ and~$\gr(p^*)$
are compatible with the realization of chord diagrams as singular
links (see Section~\ref{s:psi}).

\begin{prop}\label{p:cover} For a finite covering of connected $3$-manifolds $p:E\longrightarrow B$
the following diagram commutes. \nopagebreak $$
\begin{picture}(9.5,4.5)(-4.5,-2)
\put(-4,1){\makebox(2,1){$\Ab(E)$}}
\put(2,1){\makebox(2,1){$\gr\L(E)$}}
\put(-4,-2){\makebox(2,1){$\Ab(B)$}}
\put(2,-2){\makebox(2,1){$\gr\L(B)$}}
\put(-2,1.5){\vector(1,0){3.9}} \put(-2,-1.5){\vector(1,0){3.9}}
\put(-3,-1){\vector(0,1){2}} \put(3,-1){\vector(0,1){2}}
\put(-1,1.4){\makebox(2,1){$\psi_E$}}
\put(-1,-1.6){\makebox(2,1){$\psi_B$}}
\put(-4.8,-0.5){\makebox(1.8,1){$\pi_1(p)^*$}}
\put(3,-0.5){\makebox(1.8,1){$\gr(p^*)$}} \end{picture} $$
\end{prop}
{\bf Proof:} Let~$D$ be a~$G$-labeled chord diagram
where~$G=\pi_1(B,*)$. By relations~(Rep) we may assume that
between any two consecutive chord endpoints of~$D$ there is
exactly one label. Furthermore, we may assume that connected
components of the skeleton of~$D$ without chord endpoints have
exactly one label. Let $L_D\subset B$ be a realization of~$D$ such
that the complement of small intervals around the labels of~$D$ is
mapped into~$U$. As in the proof of Lemma~\ref{l:pstarfilt} we
write~$p^*(L_D)$ as a sum of~$n^d$ singular links~$L_i$ by local
modifications of~$p^{-1}(L_D)$ inside of~$p^{-1}(U)$. The
connected components of~$p^{-1}(U)$ are in natural one-to-one
correspondence with~$H\setminus G$. By a homotopy of singular
links we pull $L_i\cap p^{-1}(U)$ into~$\widetilde{U}$ by
following paths belonging to the lifts of elements of~${\cal R}$
with starting point~$\tilde{*}$. This way~$L_i$ becomes the
realization of an~$H$-labeled chord diagram that appears as one of
the~$n^d$ terms in the definition of~$\pi_1(p)^*(D)$.
Orient~$p^{-1}(U)$ such that~$p_{\vert p^{-1}(U)}$ is orientation
preserving. Then the path in~$E$ belonging to the lift
of~$g_a\in{\cal R}$ is orientation preserving
iff~$\sigma_G(g_a)=1$. This shows that the diagram belonging
to~$L_i$ appears in~$\pi_1(p)^*(D)$ with the correct sign and
completes the proof.
 $\Box$

\medskip

Let~$\Sigma_1$ be a ribbon graph with bands~$B_i$ ($i\in J_1$),
coupons $C_i$ ($i \in J_2$) and with a basepoint~$*$ in a
distinguished coupon~$C_o$. Let $p_1:B\longrightarrow \Sigma_1$ be
an~$I$-bundle with an orientation of~$p_1^{-1}(C_o)$. Consider a
finite connected covering~$p:E\longrightarrow B$.
Let~$\Sigma_2=p^{-1}(s(\Sigma_1))$ where~$s$ is a section
of~$p_1$. Then there exists a unique
projection~$p_2:E\longrightarrow \Sigma_2$ such that $p_1\circ
p=p_1\circ p\circ p_2$ and~$E$ is an~$I$-bundle over~$\Sigma_2$.
The surface~$\Sigma_2$ has a structure of a ribbon graph with
bands $B_i^c$ ($i\in J_1, c\in H\setminus G$) and coupons $C_i^c$
($i\in J_2, c\in H\setminus G$) such that~$p_1\circ p$ restricts
to homeomorphisms from~$C_i^c$ to~$C_i$ and from~$B_i^c$ to~$B_i$.
Given elements $a(B_1)\in\widehat{\A}(\{e\},1,+,+)$
satisfying~$a(B_1^\star)=a(B_1)^\star$ for all~$B_1\in{\cal
B}^{\rm or}(\Sigma_1)$, we define~$b(B_2)=a(p_1(p(B_2)))$ for
all~$B_2\in {\cal B}^{\rm or}(\Sigma_2)$. Choose~$\tilde{*}\in
p^{-1}(s(C_o))$. Let~$C$ be the unique coupon of~$\Sigma_2$
with~$\tilde{*}\in C$. We orient~$p_2^{-1}(C)$ such that $p_{\vert
C}$ is orientation preserving. Then we have the following theorem.

\begin{theoremX}\label{t:cover}
Let $p:E\longrightarrow B$ be a covering of
$I$-bundles~$p_1:B\longrightarrow \Sigma_1$, $p_2:E\longrightarrow
\Sigma_2$ over ribbon graphs~$\Sigma_1$, $\Sigma_2$ as above. We
assume that for all~$B_1\in{\cal B}^{\rm or}(\Sigma_1)$ we
have~$\widehat{\Delta}(a(B_1))=a(B_1)\otimes a(B_1)$. Then the
following diagram commutes.

$$
 \begin{picture}(9.5,4.5)(-4.5,-2)
\put(-4,1){\makebox(2,1){$\L(E)$}}
\put(2,1){\makebox(2,1){$\widehat{\A}(E)$}}
\put(-4,-2){\makebox(2,1){$\L(B)$}}
\put(2,-2){\makebox(2,1){$\widehat{\A}(B)$}}
\put(-2,1.5){\vector(1,0){3.9}} \put(-2,-1.5){\vector(1,0){3.9}}
\put(-3,-1){\vector(0,1){2}} \put(3,-1){\vector(0,1){2}}
\put(-1,1.4){\makebox(2,1){$Z^{p_2,b}_E$}}
\put(-1,-1.6){\makebox(2,1){$Z^{p_1,a}_B$}}
\put(-4,-0.5){\makebox(1,1){$p^*$}}
\put(3,-0.5){\makebox(2,1){$\pi_1(p)^*$}} \end{picture} $$
\end{theoremX}
{\bf Proof:} Choose~$T_1\subset \Sigma_1$ as in
equation~(\ref{e:chooseT}). Let $L$ be a generic link in~$B$. By
definition of~$Z_B^{p_1,a}$ we can express~$Z_B^{p_1,a}(L)$ in the
following form:

\begin{equation}\label{e:ZBL}
 Z_B^{p_1,a}(L)=\kappa_1\left(\bigotimes_{i\in J_2}
t_i\otimes \bigotimes_{i\in J_1} d_{s_i}\left(I(\hat{g}_i)\circ
a_{T_1}\left(B_i\right)\right)\right),
\end{equation}

where the symbols~$\kappa_1$, $t_i$, $s_i$, and $I(\hat{g}_i)$ are
explained now: $\kappa_1$ is a gluing map defined in terms of the
gluing pattern of the bands and coupons of~$\Sigma_1$. The tangles
$L\cap p_1^{-1}(C_i)$ are equipped with suitable parenthesis on
their sources and targets in the definition of~$t_i=Z(L\cap
p_1^{-1}(C_i))$. Let~$I_0^i$ be the lower base of the band~$B_i$.
Let~$C^i$ be the unique coupon of~$\Sigma_1$ with~$C^i\cap
I_0^i\not=\emptyset$. Then we define~$s_i$ as the subword
of~$\sour(L\cap p_1^{-1}(C^i))$ or~$\tar(L\cap p_1^{-1}(C^i))$
determined by~$L\cap p_1^{-1}(I_0^i)$. Finally, for
$\hat{g}_i={p_1}_*^{-1}(g_{B_i})$ we define
$I(\hat{g}_i)\in\Ab(G,\sigma_G,+,+)_0$ as a single strand without
chords labeled by~$\hat{g}_i$.

Choose~$T_2\subset\Sigma_2$ as in equation~(\ref{e:chooseT}) such
that~$(p_1\circ p)^{-1}(T_1)\subset p_2^{-1}(T_2)$. The connected
components of~$(p_1\circ p)^{-1}(T_1)$ are in a natural one-to-one
correspondence with~$H\setminus G$, so we may assume that coupons
are indexed in a way such that~$C_i^c$ lies in the part
of~$(p_1\circ p)^{-1}(T_1)$ belonging to~$c\in H\setminus G$.
Furthermore, we may assume that the bands of~$\Sigma_2$ are
indexed in a way such that the band~$B_k^c$ leads
from~$C_i^{c\cdot g_{B_k}}$ to~$C_j^{c}$ whenever the band~$B_k$
of~$\Sigma_1$ leads from~$C_i$ to~$C_j$. There exist unique
representatives~$g_c$ of cosets~$c\in H\setminus G$ such that the
lifts of the elements~$g_c$ with starting point~$\tilde{*}$ can be
represented by paths inside of~$T_2\subset E$. Using
Fact~\ref{f:assosym} we can
relate~$Z_E^{p_2,b}\left(p^{-1}(L)\right)$ to
equation~(\ref{e:ZBL}) as follows

\begin{equation}\label{e:ZE}
 Z_E^{p_2,b}\left(p^{-1}(L)\right)=\kappa_2\left(\bigotimes_{c\in
H\setminus G} \left(\bigotimes_{i\in J_2} t^c_i\otimes
\bigotimes_{i\in J_1} d_{s_i}\left(I(\hat{g}_i^c)\circ
b_{T_2}\left(B_i^c\right)\right)\right)\right),
\end{equation}

where~$\kappa_2$ is defined in terms of the gluing pattern of the
bands and coupons of~$\Sigma_2$, $t_i^c=t_i$ if~$\sigma_G(g_c)=1$
and~$t_i^c=t_i^\star$ otherwise, and
$\hat{g}_i^c={p_2}_*^{-1}(g_{B_i^c})$. We
have~$p_*(\hat{g}^c_k)=g_{c\cdot \hat{g}_k}\hat{g}_kg_c^{-1}$
because the lift of~$\hat{g}_k$ with starting
point~$\tilde{*}\cdot g_{c\cdot \hat{g}_k}$ can be represented by
a path inside of~$(p_1\circ p)^{-1}(T_1)\cup B_k^c\subset
p_2^{-1}(T_2)\cup B_k^c$ that travels through~$B_k^c$ exactly once
in the positive direction.

Let $d=[G:H]$. Define the $d$-th iterated
comultiplication~$\widehat{\Delta}_d$
by~$\widehat{\Delta}_2=\widehat{\Delta}$ and
by~$\widehat{\Delta}_d=(\widehat{\Delta}_{d-1}\widehat{\otimes}
\id)\circ \widehat{\Delta}$. Let $H\setminus G=\{c_1,\ldots,
c_d\}$. For a $G$-labeled chord diagram~$D$ we
define~$\psi_\nu(D)=D$ if~$\sigma_G(g_{c_\nu})=1$, and
$\psi_\nu(D)=D^\star$ otherwise. Using the definition
of~$\pi_1(p)^*$ with respect to the set~${\cal R}=\{g_c\;\vert\;
c\in H\setminus G\}$ we construct
$\pi_1(p)^*\left(Z_B^{p_1,a}(L)\right)$ from the elements~$t_i$
($i\in J_2$) and $d_i=d_{s_i}(I(\hat{g}_i)\circ a_{T_1}(B_i))$
($i\in J_1$) in the following three steps:

(1) for $j\in J_2$ let
$T_j=(\psi_1\widehat{\otimes}\ldots\widehat{\otimes}\psi_d)\left(\widehat{\Delta}_d(t_j)\right)$,

(2) define $D_i$ ($i\in J_1$) by replacing the labels~$\hat{g}_i$
in the $\nu$-th tensor factor of

$$(\psi_1\widehat{\otimes} \ldots\widehat{\otimes}
\psi_d)\left(\widehat{\Delta}_d(d_i)\right)$$

with $q(g_{c_\nu\cdot \hat{g}_i}\hat{g}_ig_{c_\nu}^{-1})\in H$
($\nu=1,\ldots, d$) where $q:p_*(H)\longrightarrow H$ is the
inverse of $p_*$,

(3) assume that the band~$B_k$ of~$\Sigma_1$ is leading from~$C_i$
to $C_j$. Let $\pi_k$ be the permutation determined by~$c_\nu\cdot
\hat{g}_k=c_{\pi_k(\nu)}$. Then the $\nu$-th tensor factor
of~$D_k$ is glued between the $\pi_k(\nu)$-th tensor factor
of~$T_i$ and the $\nu$-th tensor factor of~$T_j$ in the same way
as~$d_k$ is glued to~$t_i$ and to~$t_j$ by the map~$\kappa_1$.

Part~(2) of Fact~\ref{f:Zprop}
implies~$\widehat{\Delta}_d(t_j)=t_j^{\otimes d}$, so we obtain

$$T_j=t_j^{c_1}\otimes \ldots\otimes t_j^{c_d}.$$

By the assumptions on the elements~$a(B_1)$ made in the theorem we
have~$\widehat{\Delta}_d(d_i)=d_i^{\otimes d}$ and therefore

$$ D_i=\bigotimes_{\nu=1}^d d_{s_i}\left(I(q(g_{c_\nu\cdot
\hat{g}_i}\hat{g}_ig_{c_\nu}^{-1}))\circ
\psi_\nu(a_{T_1}(B_i))\right)= \bigotimes_{\nu=1}^d
d_{s_i}\left(I(\hat{g}_i^{c_\nu})\circ
b_{T_2}(B_i^{c_\nu})\right).$$

From part~(3) of the description
of~$\pi_1(p)^*\left(Z_B^{p_1,a}(L)\right)$ from above we see that
the elements~$T_j$ and~$D_i$ are sticked together according to the
gluing map~$\kappa_2$. This way we obtain the expression in
equation~(\ref{e:ZE}) for~$\pi_1(p)^*\left(Z_B^{p_1,a}(L)\right)$
which completes the proof.$\Box$

\medskip


Let~$p:S^2\times I\longrightarrow P^2\times I$ be the universal
covering of~$P^2\times I$. Recall the decomposition of~$P^2$ into
a coupon~$B_0$, a band~$B_1$ and a disk~$B$ from
Section~\ref{s:rp2}. The preimages~$p^{-1}(B_0\times 1/2
)$,~$p^{-1}(B_1\times 1/2)$, and~$p^{-1}(B\times 1/2)$ constitute
a decomposition of~$S^2$ into two coupons~$B_0^i$, two
bands~$B_1^i$ and two disks~$B^i$ ($i=0,1$). For the ribbon
graph~$S=\bigcup_{i,j=0}^1 B_i^j$ a universal Vassiliev invariant
$Z_{S\times I}=Z_{S\times I}^q$ with $q(x,t)=x$ is defined in
Section~\ref{s:vari}. Denote the inclusion map~$S\times I\subset
S^2\times I$ by~$j$. The map~$\pi_1(j)_*$ was associated
to~$\pi_1(j)$ at the end of Section~\ref{s:lcd}. It replaces all
labels of chord diagrams by the neutral element. For
links~$L\subset S\times I\subset S^2\times I$ we
define~$Z_{S^2\times I}(L)=\pi_1(j)_*\left(Z_{S\times
I}(L)\right)$. Any link in~$S^2\times I$ is isotopic to a link
in~$S\times I$. As a corollary to part~(1) of the following
theorem we see that the definition of~$Z_{S^2\times I}$ leads to
an isotopy invariant of links in~$S^2\times I$.

\begin{theoremX}\label{t:cover2}
(1) Let $i:S^2\times I\longrightarrow \R^3$ be the inclusion
defined by $i(x,t)=(t+1)x$. Then for links~$L\subset S\times
I\subset S^2\times I$ we have

$$Z_{S^2\times I}(L)=Z(i(L)).$$


(2) For the universal covering $p:S^2\times I\longrightarrow
P^2\times I$ we have

$$\pi_1(p)^*\circ Z_{P^2\times I}=Z_{S^2\times I}\circ p^*.$$
\end{theoremX}
{\bf Proof:} (1) Let $L\subset S\times I$ be a link in standard
position. The inclusion~$i_{\vert S\times I}:S\times
I\longrightarrow \R^3$ is depicted on the left side of
Figure~\ref{f:ribbonS}.

\begin{figure}[!h] \centering
\setbox1=\hbox{\input{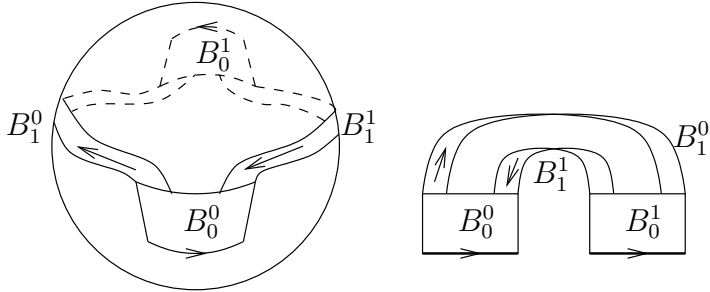}}
$\vcenter{\box1}$ \caption{The ribbon graph~$S\subset S^2$ and the
abstract ribbon graph~$S$} \label{f:ribbonS}
\end{figure}

We can regard~$i(L\cap B_\mu^\nu\times I)$ ($\nu, \mu=0,1$) as
non-associative tangles such that

$$ Z(i(L))=Z(i(L\cap B_0^1\times I))\circ \left(Z(i(L\cap
B_1^0\times I))\btimes Z(i(L\cap B_1^1\times I))\right)\circ
Z(i(L\cap B_0^0\times I)),$$

the non-associative tangle $i(L\cap B_1^0\times I)$ (resp.\
$i(L\cap B_1^1\times I)$) consists of a left-handed (resp.\
right-handed) half twist of a bunch of strands, and we have

$$ s_\nu=\sour(i(L\cap B_1^\nu\times I))=\tar(i(L\cap
B_1^\nu\times I))^\vert \quad\mbox{($\nu=0,1$).}$$

We compute the invariants~$Z(i(L\cap B_1^\nu\times I))$ using the
analytic definition of~$Z$: from Example~1.4, Property~1.12
with~$A=-1$, and Section~3.3 of~\cite{Les} we obtain the left side
of equation~(\ref{e:Zhalftwist}) for some~$u\in\widehat{\A}(\{e\},
1, s_\nu, s_\nu)$.

\begin{equation}\label{e:Zhalftwist}
Z(i(L\cap B_1^\nu\times I))\quad =\quad \mbox{\setbox1=\hbox{\input{ht1}}
$\vcenter{\box1}$}\quad
=\quad \mbox{\setbox1=\hbox{\input{ht2}}
$\vcenter{\box1}$}
\end{equation}

where $t_\nu=-(-1)^\nu\mbox{\setbox1=\hbox{\begin{picture}(0,0)%
\includegraphics{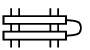}%
\end{picture}%
\setlength{\unitlength}{3947sp}%
\begingroup\makeatletter\ifx\SetFigFont\undefined%
\gdef\SetFigFont#1#2#3#4#5{%
  \reset@font\fontsize{#1}{#2pt}%
  \fontfamily{#3}\fontseries{#4}\fontshape{#5}%
  \selectfont}%
\fi\endgroup%
\begin{picture}(396,211)(1039,-410)
\end{picture}
}
$\vcenter{\box1}$}/4$. The right side of
equation~(\ref{e:Zhalftwist}) follows because we can use the
relations in Figure~\ref{rels} to cancel~$u$
with~${u^{\vert}}^{-1}$. Let $b_1^\nu$ be the part of degree~$0$
of~$Z(i(L\cap B_1^\nu\times I))$. Again by applying the relations
in Figure~\ref{rels} we obtain

\begin{equation}\label{e:ZiL}
Z(i(L))=Z(i(L\cap B_0^1\times I))\circ \left(b_1^0\btimes
b_1^1\right)\circ Z(i(L\cap B_0^0\times I)). \end{equation}

Recall that in the definition of~$Z_{S^2\times I}(L)$ the
invariants of non-associative tangles~$Z(L\cap B_0^\nu\times I)$
are sticked together using a gluing map determined by a diagram of
the ribbon graph~$S$ shown on the right side of
Figure~\ref{f:ribbonS}. This implies

\begin{equation}\label{e:ZS2L}
Z_{S^2\times I}(L)=Z(L\cap B_0^1\times I)^{-\vert}\circ
\left(b_1^0\btimes b_1^1\right)\circ Z(L\cap B_0^0\times I),
\end{equation}

where the symmetry operations~$-$ and~$\vert$ have been defined at
the end of Section~\ref{s:Z}. Since~$i(L\cap B_0^1\times I)=(L\cap
B_0^1\times I)^{-\vert}$ and~$i(L\cap B_0^0\times I)=L\cap
B_0^0\times I$ we obtain the first part of the theorem from
Fact~\ref{f:assosym} and equations~(\ref{e:ZiL})
and~(\ref{e:ZS2L}).


\smallskip

(2) Let $X=B_0\cup B_1\subset P^2$. Let~$q=p_{\vert S\times
I}:S\times I\longrightarrow X\times I$. For a link~$L$ in~$X\times
I$ Theorem~\ref{t:cover} implies

\begin{equation}\label{e:pZold}
\left(\pi_1(q)^*\circ Z_{X\times I}\right)(L)=Z_{S\times
I}\left(p^{-1}(L)\right).
\end{equation}

Let $\pi_1(j)_*$ be the map from above given by forgetting labels
of a chord diagram. Denote the inclusion map~$X\times I\subset
P^2\times I$ by~$k$. Then~$\pi_1(k)_*$ is given by reduction
modulo~$2$ of the labels of a~$\Z$-labeled chord diagram, where
$\Z\cong\pi_1(X\times I,*)$. Since it is easy to verify

$$\pi_1(j)_*\circ\pi_1(q)^*=\pi_1(p)^*\circ\pi_1(k)_*$$

we obtain part~(2) of the theorem by applying $\pi_1(j)_*$ to both
sides of equation~(\ref{e:pZold}).
 $\Box$

\medskip

It is easy to show that the definition of~$Z_{S^1\times S^2}$
(resp.\ $Z_{S^1\tilde{\times} S^2}$) in Section~\ref{s:s1s2}
extends to the case where~$Z_{S^1\times I\times I}$ (resp.\
$Z_{X\times I}$) is defined using an arbitrary ribbon graph
structure on~$S^1\times I$ (resp.\ $X\times I$). Let $B=S^1\times
S^2$ or $B=S^1\tilde{\times}S^2$. For a finite covering
$p:E\longrightarrow B$ with a connected space~$E$ we have~$E\cong
S^1\times S^2$ if $B=S^1\times S^2$ and $E\cong S^1\times S^2$ or
$E\cong S^1\tilde{\times}S^2$ otherwise. The map~$\pi_1(p)^*$
descends to a map~${\cal E}(B)\longrightarrow{\cal E}(E)$ and a
statement similar to Theorem~\ref{t:cover} holds for finite
coverings of~$B$.


{\small\tt Jens Lieberum, Mathematisches Institut, Rheinsprung 21,
CH-4051 Basel,

email:\ lieberum@math.unibas.ch}

}

\end{document}